\documentstyle{amsppt}
\NoBlackBoxes
\def\Det{\operatorname{Det}}
\def\D{\operatorname{Det}}
\def\la{\lambda}
\def\C{\Bbb C}
\def\Z{\Bbb Z}

\def\ga{\gamma}
\def\De{\Delta}
\def\Dif{\operatorname{Dif}}
\def\Lie{\operatorname{Lie}\,}
\def\SL{\operatorname{SL}}

\def\a{\alpha}

\def\fn{{\frak n}}
\def\fh{{\frak h}}
\def\ad{\operatorname{ad}}
\def\I{Ind_\mu}
\def\In{Ind_a}
\def\g{gl}
\def\e{\operatorname{deg}}
\def\n{\frak n}
\def\A{\frak A}
\def\lnorm{||}
\def\rnorm{||}
\def\lmod{\left |}
\def\rmod{\right|}
\def\BI{\Bbb I}
\def\1{\operatorname{I}\!\!\!\!1}

\topmatter

\title
Certain Topics on the Lie Algebra
\lowercase{$gl(\lambda)$} Representation\\
Theory
\endtitle

\address\endaddress
\email borya\@mccme.ru\qquad
borya\@main.mccme.rssi.ru\qquad
shoykhet\@ihes.fr
          \endemail


\author
B.~B.~Shoikhet
\endauthor


\endtopmatter

\document

\subhead
Introduction
\endsubhead
\subsubhead
0.1
\endsubsubhead
The aim of this work is to give a review of recent results in the
representation theory of the Lie algebra $gl(\la)$ and their applications
to the new combinatorial identities (see 0.11).

The $gl(\la)$ is an infinite-dimensional Lie algebra,
dependent on the parameter $\la\in \C$, which is a continuous version
of the Lie algebra $gl_\infty$. It was introduced by B\.L.~Feigin in~[1]
for calculating the homology ring of the Lie algebra of differential operators
on a line and has several equivalent definitions. Here we shall
consider three of them.

\subsubhead
0.2
\endsubsubhead
First, it is the Lie algebra constructed from the associative algebra
of twisted differential operators (see~[7]) on $\C P^1$. The differential
operator of order $\le k$ on the line bundle $\Cal L$ is a global object
which locally acts on the holomorphic sections of this bundle and
$[\ldots[[\Cal D,f_1],f_2]\ldots f_{k+1}]=0$ for any $k+1$ holomorphic
functions (these functions are considered as zero-order operators
$f_i:\Gamma(\Cal L)\to\Gamma(\Cal L)$). Any line bundle on $\C P^1$
is $\Cal O(n)$ for a certain $n\in\Z$, and we obtain for each $n\in\Z$
an associative algebra of twisted differential operators (for
$n=0$ they are ordinary differential operators). In fact, such an algebra
of twisted differential operators exists for any $n=\la\in \C$ (irrespective
of the fact that the corresponding line bundle does not exist for a nonintegral
$\la$).

Recall the construction of this algebra (see [6, 7] for details). We denote
by  $D_{\Cal O(1)}$ a sheaf of (ordinary) differential
operators on the total space of the bundle $\Cal O(1)$ on $\C P^1$, and
let $E$ be a  Euler vector field on fibers of this bundle. We denote by $D_E$
the subsheaf in $D_{\Cal O(1)}$ consisting of operators commuting with $E$,
and then
$\{\Cal D E-\la\Cal D\mid\Cal D\in D_E\}$
is a two-sided ideal in
$D_E$, and we denote by $\Dif_\la$ the (associative) quotient algebra of
$D_E$ over this ideal. Then for $\la\in\Z$ $Dif_\la$
coincides with the algebra of
differential operators on the bundle $\Cal O(\la)$. We can  give the
following heuristic explanation for this fact. $Dif_\la$ are differential
operators which must act on the sections of the ``degree of homogeneity
$\la$'' $(\la\in\C)$.
Using two holomorphic sections $\ga_1$ and $\ga_2$ of the
bundle $\Cal O(1)$, we can construct two functions
$\tilde{\ga}_1$ and $\tilde{\ga}_2$, linear
along the fibers, on which $E$ acts  with the eigenvalue~1.
>From the sheaf
$\tilde{\ga}_1^\la\cdot{\Cal O}_{\Cal O(1)}+
\tilde{\ga}_2^\la\cdot{\Cal O}_{\Cal O(1)}$ we must take all sections of
the degree of homogenetuty $\la$, and $Dif_\la$ is what acts on
these sections. The subalgebra $D$ consists of differential operators which
preserve the proper decomposition of the operator $E$, and the
factorization with respect to the ideal
$\{\Cal D E-\la\Cal D\mid\Cal D\in D_E\}$  is the isolation of the component
corresponding to the eigenvalue $\la$.

>From this point of view $gl(\la)=\Lie(Dif_\la)$.

\subsubhead
0.3
\endsubsubhead
We have the representation
$p_\la:U(sl_2)\to Dif_\la(\C P^1)$
for all
$\la\in\C$ (see [6, 7]). In order to construct it, we represent $\Cal O(1)$
as $\SL_2(\C)\underset{B}\to{\tilde\times}{\widetilde\C}$, where
$\widetilde{\C}$ is the
corresponding 1-dimensional representation of the Borel subgroup
$B\subset \SL_2$. Then $\SL_2(\C)$ acts on the total space of the bundle
$\Cal O(1)$, and the Lie algebra $sl_2(\C)$ is mapped into the vector fields
on $\Cal O(1)$. This is a Lie algebra homomorphism, and therefore, the map
$p_\la:U(sl_2)\to\Dif_\la$  is well-defined.
The Beilinson--Bernstein theorem~[6]
states in this simplest case that $p_\la$ is surjective and that the kernel of
$p_\la$ is a two-sided ideal in $U(sl_2)$ generated by
$\Delta-\frac{\la(\la+2)}2$,
where
$\Delta=e\cdot f+f\cdot e+\frac{h\cdot h}{2}\in U(sl_2)$ is the Casimir
operator~[4].  On the Verma $sl_2$-module with the highest weight
$\la$ the operator $\Delta$
acts  by the multiplication by
$\frac{\la(\la+2)}2$.
Thus
$gl(\la)=\Lie\Bigl(U(sl_2)\Big/\Bigl(
\frac{\la(\la+2)}2\Bigr)\Bigr)$.

\subsubhead
0.4
\endsubsubhead
The name of the Lie algebra $gl(\la)$ itself implies a connection with the
Lie algebra $gl_n$. We shall now give its last definition from which it is
clear, in particular,  in what sence $gl(\la)$ is a Lie algebra of matrices
of complex dimension.

Let us consider an $n$-dimensional irreducible $sl_2$-module $V$. There are
a surjection
$U(sl_2)\to gl(V)$
and a principal $sl_2$-subalgebra
$sl_2\subset U(sl_2)\to gl(V)$ in $gl(V)$. Being a $sl_2$-module,
$gl(V)\cong gl_n\cong\oplus_{i=0}^{n-1}\pi_i$, where $\pi_i$ is an
irreducible $(2i+1)$-dimensional $sl_2$-module. Being a Lie algebra, $gl_n$
is generated by $\pi_0$, $\pi_1$, and $\pi_2$, and for a sufficiently large
$n$ relations of a fixed degree, depend on $n$ analytically (these relations
were described in~[1]). This makes it possible to consider $n$ to be a
complex number $\la$ by substituting $\la$ for $n$ in all relations, and
 in this way define the structure of the Lie algebra on
$\oplus_{i=0}^{\infty}\pi_i$ dependent on $\la\in \C$.
It is easy to show (see~[1] and 1.1 in Ch.~2) that this definition
is equivalent to the preceding ones.

\subsubhead
0.5
\endsubsubhead
The Lie algebra $gl(\la)$ is $\Z$-graded, i.e.,
$gl(\la)^i=\{\xi\in gl(\la)\mid[h,\xi]=2i\xi\}$, where
$h\in sl_2\subset gl(\la)$, and there exists the Cartan decomposition
(see~[14]) $gl(\la)={\frak n}_-\oplus{\frak h}\oplus{\frak n}_+$, where
${\frak n}_-=\oplus_{i<0}gl(\la)^i$, $\frak h=gl(\la)^0$,
${\frak n}_+=\oplus_{i>0}gl(\la)^i$. This Cartan decomposition is
consistent with the decomposition of $gl_n$ (see 0.4). However, an
essential difference from the classical Lie algebra is that ${\frak n}_{\pm}$
cannont be decomposed into a direct sum of one-dimensional $h$-invariant
subspaces. For instance, ${\frak n}_+^{(1)}$ is entirely $\frak h$-invariant
(it cannont be decomposed even into a direct sum of two proper subspaces)
and is a space of simple positive root vectors.  From this point of view
$gl(\la)$ is a Lie algebra with a continual system of roots.

This can be illustrated by a simple example. A parabolic subalgebra in
$gl(\la)$ corresponding to the roots $\a_1,\ldots,\a_k\in\C$ is the
subalgebra ${\frak p}_{\a_1\ldots,\a_k}$ generated by
${\frak n}_+$, $\frak h$, and
${\tilde{\frak n}}_-^{(-1)}=\{P(h)f\mid P(h)=
(h-a_1)\cdots$ $(h-\a_k)P_1(h),\,P_1\in\C[h]\}$. (Here we use the second
definition of $gl(\la)$ as
$\Lie\Bigl(U(sl_2)\Big/\Bigl(\De-
\frac{\la(\la+2)}2\Bigr)\Bigr)$, see 0.3.) For $\la\in\Z$,
there exists a surjection $gl(\la)\to gl_{\la+1}$ (see 0.4) and under this
surjection almost
all these parabolic subalgebras pass into one (independent of
$\a_1,\ldots,\a_k$ ) parabolic subalgebra in $gl_{n+1}$. In the
classical case, parabolic subalgebras correspond to subsets of simple
positive roots.

\subsubhead
0.6
\endsubsubhead
The subject of our study is representations with the highest weight
(see~[8, 14, 4]) of the Lie algebra $gl(\la)$ all of whose levels
(relative to the $\Z$-grading, see 0.5) are finite-dimensional. Thus, we
study the representations $V$ of the Lie algebra $gl(\la)$ which satisfy the
following conditions:

(1) there exists a vector $v\in V$ such that
$$
{\frak n}_+v=0,\ \text{ and }\ hv=\chi(h)v\ \text{ for all }\ h\in\frak h
$$

(2) $V=U(gl(\la))\cdot v$

(3) all levels of $V$ are finite-dimensional.

Clearly, this representation  exists not for all $\chi\in{\frak h}^*$.
For the general $\chi$, the corresponding Verma module is irreducible
(and has infinite-dimen\-sio\-nal levels). In fact, any representation of this
kind with the $k$-dimen\-sio\-nal first level is a quotient of the representation
induced with a certain character from ${\frak p}_{\a_1,\ldots,\a_2}$
(see 1.3 in Ch.~2). The space of characters is $k$-dimensional and we
obtain a $(2k)$-parametric family of representations. For the general
parameters, these representations are irreducible, and our first objective
is to find the parameters for which these representations are reducible.

In this review, we only consider a 2-parameter family of representations
with 1-dimensional first level.

\subsubhead
0.7
\endsubsubhead
We follow the method of Kac and Radul~[2] who studied the representations
of a Lie algebra close to $gl(\la)$, namely that of the
differential operators on a circle.

The inclusion $\varphi_s:gl(\la)\hookrightarrow gl_{\infty,s}$, where
$gl_{\infty, s}$ is a Lie algebra analytically depending on $s\in\C$ and
isomorphic for the general $s$ to the Lie algebra $gl_\infty$ of generalized
Jacobian matrices, is defined for any $s\in\C$ (see 2.1 of Ch.~2).
It turns out that $\varphi_s$ can be extended to the mapping
$\varphi_s^{\Cal O}:gl^{\Cal O}(\la)\to gl_{\infty,s}$, which is surjective.
Here the Lie algebra $gl^{\Cal O}(\la)$ is a certain completion
of $gl(\la)$, and the $gl(\la)$-invariant subspace of any
$gl_{\infty,s}$-module is also $gl^{\Cal O}(\la)$-invariant (and hence
$gl_{\infty,s}$-invariant).

This reduces the problem of the irreducibility of $gl(\la)$-modules to the
corresponding problem of the irreducibility of
$\widehat{gl}_{\infty,s}$-modules,
and the latter problem can easily be reduced to the problem for
$\widehat{gl}_\infty$. (We easily pass to central extension since
$H^2(gl(\la),\C)=0)$.

\subsubhead
0.8
\endsubsubhead
In Ch.~1, we consider a ``model'' situation, namely, the representation
$Ind_\mu$  of the Lie algebra $\widehat{gl}_\infty$ induced from the largest
parabolic subalgebra. This representation has a zero highest weight and a
central charge $\mu\in\C$. For the modules  $Ind_\mu$, we find all
singular vectors and describe the Jantzen filtration (see~[13]) in terms
of irreducible
$gl_{\frac{\infty}{2}}^{(1)}\oplus gl_{\frac{\infty}{2}}^{(2)}$-modules.
(Here $gl_{\frac{\infty}{2}}$ and $gl_{\frac{\infty}{2}}$  are the
``upper'' and ``lower'' subalgebras in $\widetilde{gl}_\infty$,
see Fig.~1 in Ch.~1.)
We also find the determinant of Shapovalov's form (see~[15]) on all levels as
a function of $\mu$. As an obvious consequence, we obtain formulas for the
character of the irreducible first consecutive quotient module
of the Jantzen filtration.
For example, in the case of $\mu=1$ we get the classical Euler identity
(see~[10])
$$
\frac{1}{\prod_{i\ge 1}(1-q^i)}=1+\sum_{k\ge 1}\frac{q^{k^2}}
{(1-q)^2\cdot\,\ldots\,\cdot(1-q^k)^2};
$$
in the case of $\mu\in\Z_{>1}$, we get the ``higher'' Euler identities and,
in the case of $\mu\in\Z_{\le -1}$, we get expressions for the characters of
the corresponding irreducible representations. Thus, for $\mu=-1$, we find
that the character of the irreducible $\widehat{gl}_\infty$-module
with zero highest
weight and the central charge $-1$ is
$$
\chi_{-1}=1+\sum_{k\ge 1}\frac{q^k}{(1-q)^2\cdot\,\ldots\,\cdot(1-q^k)^2}.
$$

We prove that the consecutive quotient modules of the Jantzen filtration of the
representation $Ind_\mu$ are simple and that the terms of the Jantzen
filtration exhaust all submodules of $Ind_\mu$; in particular, we obtain
explicit expression for the characters of the highest irreducible consecutive
quotient modules of the Jantzen filtration.

Considering the Lie algebra $gl_{2n}$
instead of $\widehat{gl}_\infty$, we get the
``finite forms'' of all the identities and formulas.

Another expression for the determinant of Shapovalov's form of the
representation $Ind_\mu$ was obtained by Jantzen (see~[13]).

\subsubhead
0.9
\endsubsubhead
There exist $\widehat{gl}_{\infty,s}$-modules $Ind_{\mu,s}$ similar to the
$\widehat{gl}_\infty$-modules  $Ind_\mu$. Considering the inverse images
$\theta_s^*(Ind_{\mu,s})$ of the $\widehat{gl}_{\infty,s}$-modules
$Ind_{\mu,s}$ under the embedding
$\theta_s:gl(\la)\hookrightarrow\widehat{gl}_{\infty,s}$, we obtain a 2-parameter
family of representations of the Lie algebra $gl(\la)$ with a 1-dimensional
first level. As was noted in 0.7, we can solve the problem of irreducibility
for these representations. It turns out, however, that in this way we can get
representations induces from all parabolic subalgebras of the corresponding
dimension, except for two (for the general $\la$). The determinant of
Shapovalov's form of the representation $\theta_s^*(Ind_{\mu,s})$ is a
product, and when a parabolic subalgebra approaches the exceptional one,
some factors have a simple zero or pole. It follows from the
irreducibility theorem (for the general $\la$)
of the representation of $gl(\la)$
induced from two exceptional parabolic subalgebras as that the order of the
zero is equal to the order of the pole. Equating the corresponding numbers
on all levels, we get a local identity (see 0.11).

Next, in order to prove the irreducibility theorem  of the
representations of $gl(\la)$ induced from exceptional subalgebras (for the
general $\la$), we consider the embedding
$\theta_s:gl(\la)\hookrightarrow\widehat{gl}_{\infty,s}(\C[t]/t^2)$; the
highest weight of the representation induced from an exceptional paraaboloic
subalgebra is equal to the highest weight of the inverse image, under
$\theta_s$, of a certain induced representation
$I$ of $\widehat {gl}_{\infty,s}(\C[t]/t^2)$, and since the completions
$\theta_s^{\Cal O}:gl(\la)\hookrightarrow\widehat{gl}_{\infty,s}(\C[t]/t^2)$
are surjective, it follows, by analogy with Subsec~0.7, that this reduces the
problem to the calculation of the charaacter of the irreducible
quotient of the module $I$. This is the subject of Ch.~2.
\subsubhead
0.10
\endsubsubhead
In Ch.~3, we consider the Lie algebras $gl(\la)$, $\la\in\C$, as those
corresponding to the finite points of the Riemannian sphere $S^2$. In  this
case, in the neighborhood of the point $\{\infty\}\in S^2$, the Lie algebra
$gl(\la)$ is deformed into the Lie algebra of regular functions on a
nondegenerate symplectic leaf of a standard foliation in $sl_2^*$ with an
induced Poisson bracket (Lie algebras of functions for all nondegenerate
leaves are isomorphic). Thus, we assume that the point $\{\infty\}\in S^2$
is associated with a Lie algebra and the whole family of Lie algebras on
$S^2$ is decomposed into an infinite direct sum of line bundles on $S^2$.

By choosing, at every point, an induced representation of the corresponding
Lie algebra, which depends holomorphically on the point, we can obtain a
situation in which  an arbitrary induced representation of the Lie
algebra of functions on a hyperboloid sits at the point
$\{\infty\}$. Then the
determinant of Shapovalov's form of some level is a holomorphic section
of a certain line bundle on $S^2$, its Chern class can be easily found.
On the other hand, this Chern class is equal to the sum of zeros with
multiplicities of the determinant of Shapovalov's
form over all points of the sphere (on this level). We prove the
irreducibility theorem
of induced representations of the Lie algebra of
functions on a hyperboloid for the generic values of the parameters, and then
the sum is extended only to the finite points of the sphere at which we can
calculate these multiplicities by the methods given in Chs.~1 and 2.
Combining these calculations for all levels, we obtain the {\it global
identity\/} (see Subsec~0.11).

In fact, we prove the irreducibility of induced representations (for the
general parameters) of the Lie algebra of functions on a cone, that is,
on a degenerate leaf of a foliation in $sl_2^*$,
and this implies an assertion  for nondegenerate leaves
(the Lie algebras of functions on all nondegenerate
leaves are isomorphic). The Lie algebra of functions on a cone is
nilpotent up to with  subalgebra $sl_2$, i.e., the
functions that have the point $0\in\C^3$  a zero of a certain degree
form an ideal at it, and we apply Kirillov's theory~[5] according to which
representations of nilpotent Lie algebras induced from the largest
subalgebras are irreducible.

\subsubhead
0.11
\endsubsubhead
The local identity:
$$
\frac{d}{da}\biggl[\prod_{i\ge 1}\frac1{(1-a\cdot q^i)^i}\biggr]
\bigg|_{a=1}=\sum\Sb{\text{over all}}\\{\text{Young diagrams}}\\
{\Cal D}_{l_1,\ldots,l_k}\endSb \#{\Cal D}_{l_1,\ldots,l_k}\cdot
q^{\sum l_i\cdot i^2}\cdot\bigl(\chi({\Cal D}_{l_1,\ldots,l_k})\bigr)^2.
$$

The global identity:
$$
\gather
\frac{d}{da}\biggl[\frac1{(1-q)}\cdot\frac1{(1-aq^2)(1-a^2q^2)}\cdot
\frac1{(1-a^2q^3)(1-a^3q^3)(1-a^4q^3)}\cdot\,\dots\biggr]\bigg|_{a=1}=\\
\frac{d}{da}\biggl[\prod_{i\ge 1}\frac1{(1-a\cdot q^i)^i}\biggr]
\bigg|_{a=1}+\\
2\sum_{k_+\ge 1}\frac{d}{da}\biggl[\prod_{i=1}^{k_+}\frac1{(1-q^i)^i}\cdot
\prod_{i=k_++1}^{\infty}\frac{1}{(1-q^i)^{k_+}(1-aq^i)^{i-k_+}}\biggr]
\bigg|_{a=1}-\\
\sum\Sb{\text{over all Young}}\\{\text{diagrams}\
{\Cal D}_{l_1,\ldots,l_k}}\endSb\left\{\foldedtext\foldedwidth{1.5in}
{the length of the ``central'' diagonal  ${\Cal D}_{l_1,\ldots,l_k}$}
\right\}\cdot q^{\sum l_i\cdot i^2}
\cdot\bigl(\chi({\Cal D}_{l_1,\dots,l_k})\bigr)^2.
\endgather
$$

Here ${\Cal D}_{l_1,\ldots,l_k}$ is the Young diagram consisting of blocks
$1\times l_1, 2\times l_2,\ldots,k\times l_k$; the ``central''
diagonal is the diagonal beginning at its upper left vertex (see Fig.~1),
$\#{\Cal D}_{l_1,\ldots,l_k}$ is the number of squares in
${\Cal D}_{l_1,\ldots,l_k}$.
\midinsert
\vspace{+5cm}
\centerline
{Fig.~1. The Young diagram ${\Cal D}_{l_1,\ldots,l_k}$ and its ``central''
diagonal}
\endinsert

$\chi({\Cal D}_{l_1,\ldots,l_k})$
is the corresponding ``semiinfinite''
character defined as follows: consider the Lie algebra
$gl_{\frac{\infty}2}$
of the matrices $(a_{ij})$, $i$, $j=1,\ldots,\infty$ and let
$\a_1^{\lor},\a_2^{\lor},\ldots$ be its simple coroots. Then
$\chi({\Cal D}_{l_1,\ldots,l_k})$ is the character of the irreducible
$gl_{\frac{\infty}2}$-module with the highest weight $\chi$ such that
$\chi(\a_1^{\lor})=l_1,\ldots$, $\chi(\a_k^{\lor})=l_k$,
$\chi(\a_{k+1}^{\lor})=\ldots=0$. For example, if $\Cal D$ is a single block
$1\times k$ or $k\times 1$, then $\chi(\Cal D)=
\frac1{(1-q)\cdot\,\ldots,\,(1-q^k)}$. In general,
$\botsmash\chi({\Cal D}_{l_1,\ldots,l_k})$
can be easily found from Weyl's formula~[4].

The left-hand sides of these identities appear in combinatorial identities
connected with plane partitionings (see~[10]]).
\subsubhead
0.12
\endsubsubhead
This work is based on numerous talks I had with B.~L.~Feigin;
without them, it could not come into existence.
I express my deepest gratitude to him.

I am also grateful to V.\,A.~Ginzburg for a discussion of Jantzen's
filtration, and to S.\,M.~Arkhipov for his interest in my work.

I express my gratitude to ISF for partial financial support
through the grant N9R000.

\head
Chapter~1.\\
Representations of the Lie Algebras \lowercase{$\g_{2n}$} and
\lowercase{$\widehat{gl}_\infty$}\\
Induced from the Largest Parabolic Subalgebra
\endhead
\subhead
1. Reducibility and Singular Vectors
\endsubhead
\subsubhead
1.0
\endsubsubhead
There are two commuting subalgebras $gl_{n}$ in the Lie algebra $gl_{2n}$,
namely, the ``upper'' subalgebra and the ``lower'' one (see Fig.~1). We
denote them by $gl_n^{(1)}$ and $gl_{n}^{(2)}$. There are also two Abelian
subalgebras ${\frak A}_+$ and ${\frak A}_-$.
\midinsert
\vspace{+5cm}
\centerline{Fig.~1. }
\endinsert

Let $\a_{-n+1}^\lor,\ldots,\a_{-1}^\lor;\a_0^\lor;
\a_1^\lor,\ldots,\a_{n-1}^\lor$ be simple coroots. We shall consider the
character $\chi_\mu:\frak h\to\C$ $(gl_{2n}={\frak n}_-
\oplus\frak h\oplus{\frak n}_+)$ such that  $\chi_\mu(\a_0^\lor)=\mu$,
$\chi_\mu(\a_i^\lor)=0$ for $i\ne 0$. It is clear that this is a
necessary and sufficient condition for determining the representation
$Ind_\mu$ induced from $\frak p=gl_n\oplus{\frak A}_+\oplus gl_n$ with
the highest weight $\chi_\mu$. Thus, we have defined the family of
representations $Ind  _\mu$ $(\mu\in\C)$ of the Lie algebra $gl_{2n}$.
They are the main object of our study in Ch~1. We introduce certain notations:
$x_{ij}$ and $y_{ij}$ are elements of ${\frak A}_+$ and ${\frak A}_-$,
respectively $(i$, $j=1\ldots n)$, in the case of ${\frak A}_+$,
the numbering goes upward and to the right and, in the case of ${\frak A}_-$,
it goes downward and to the left. For example, $y_{11}$ is the only element
of ${\frak A}_-$ graded $-1$.

Next, let $A_k$ be a $k\times k$ matrix in the upper right corner in
$A_k=(y_{ij})$, $i,\,j=1\ldots k$, and $\Det_k\in{\Cal U}({\frak A}_-)\cong
S^*({\frak A}_-)$  be the determinant of the matrix $A_k$. For example,
$\Det_1=y_{11}$, $\Det_2=y_{12}\cdot y_{21}-y_{11}\cdot y_{22}$ and so on.
Let $v$ be the highest weight vector of $Ind_\mu$.

{\bf 1.1.} The main result of Sec.~1 is the following theorem
\proclaim
{Theorem~1}
The representations $Ind  _\mu$ of the Lie algebra $gl_{2n}$ are reducible
for $\mu\in\Z_{\ge-n+1}$ and irreducible for other $\mu\in \C$. For
$\mu\in\Z_{\ge 0}$, the singular vectors in $Ind  _\mu$ are
$$
\Det_1^{\mu+1}\cdot v,\Det_2^{\mu+2}\cdot v,\ldots,\Det_n^{\mu+n}\cdot v.
$$
For $-n+1\le\mu\le 0$ and $\mu\in\Z$, the singular vectors in $Ind  _\mu$ are
$$
\Det_{-\mu+1}\cdot v,\Det_{-\mu+2}^2\cdot v,\ldots,
\Det_n^{\mu+n}\cdot v.
$$
\endproclaim
The remaining part of Sec.~1 is devoted to the proof of this theorem.
\remark
{Remark}
$\Det_k$ are weight elements of $\Cal U({\frak A}_-)$ of the weight,
$-\bigl(k\a_0+(k-1)(\a_{-1}+\a_1)+\dots(\a_{-k+1}+\a_{k-1})\bigr)$. For the
$\1$-specialization, this weight is equal to $-k^2$.
\endremark
\subsubhead
1.2
\endsubsubhead
Let $gl_n^{(1)}=\fn_-^{(1)}\oplus\fh^{(1)}\oplus\fn_+^{(1)}$,
$gl_n^{(2)}=\fn_-^{(2)}\oplus\fh^{(2)}\oplus\fn_+^{(2)}$, and $v$ be the
highest weight vector in $Ind  _\mu$. Any element of $Ind  _\mu$ has the form
$\xi\cdot v$, where $\xi\in\Cal U({\frak A}_-)=S^*({\frak A}_-)$. If
$\xi\cdot v$ is a singular vector, then
$$
[e,\xi]\cdot v\in\Cal U({\frak A}_-)\cdot(\fn_-^{(1)}\oplus\fn_-^{(2)})v,
$$
where $e\in\fn_+\subset gl_{2n}$ is arbitrary.

The main observation is that $\Cal U(\frak A)$ is {\it invariant\/} with
respect to the $\ad$-action of $\fn_+^{(1)}$ and $\fn_+^{(2)}$ i.e.,
we searh for vectors in
$Ind_\mu$ which are singular for $gl_n^{(1)}\oplus gl_n^{(2)}$ and then find
the values of $\mu$ when these singular vectors pass to 0 uppon the application
of $e_0=x_{11}$. To be more precise, if $e_{-1},\ldots,e_{-n+1}$
$(e_1,\ldots, e_{n-1})$ are root vectors from $gl_n^{(1)}$ $(gl_n^{(2)})$
corresponding to the simple positive roots
$\a_{-1},\ldots,\a_{-n+1}(\a_1,\ldots,\a_{n-1})$, then the operators
$\ad(e_i)$ for $i<0$ shift the $(|i|+1)$th column of ${\frak A}_-$ (reckoning
from the center) to the right by 1 and for $i>0$ they shift the $(i+1)$th row
upward by 1. In other words, on the space
$\Cal U({\frak A}_-)=\C[y_{ij};i,j=1,\ldots,n]$ $\ad(e_i)$  for $i<0$
acts by the vector field
$\sum\limits_{s=1\ldots n}y_{s,|i|}\frac{\partial}{\partial y_{s,|i|+1}}$ and for
$i>0$ by the vector field $\sum y_{i,s}\frac{\partial}{\partial y_{i+1,s}}$.
\subsubhead
1.3
\endsubsubhead
The following lemma refers to the action of the operators $\ad(\fn_+^{(1)})$
and $\ad(\fn_+^{(2)})$ on $\Cal U({\frak A}_-)$.
\proclaim
{Lemma~1}
The monomials $\Det_1^{k_1}\cdot\,\ldots,\cdot\Det_n^{k_n}\cdot v$ and only
these vectors are
$gl_n^{(1)}\oplus gl_n^{(2)}$-singular vectors in
$Ind  _\mu$ \rom(for any $\mu$\rom).
\endproclaim
\demo
{Proof}
The fact that these monomials are
$gl_n^{(1)}\oplus gl_n^{(2)}$-singular vectors follows from the results of
Subsec.~1.2. Conversely, let $\xi\cdot v$  be a
$gl_n^{(1)}\oplus gl_n^{(2)}$-singular vectors in
$Ind  _\mu$ ($\xi\in S^*({\frak A}_-)$). We shall consider the {\it minimal\/}
rectangular domain in ${\frak A}_-$  with the vertex in the central (upper
right) corner which includes all elements which are contained in the
notation of $\xi$. Suppose, for example, that it is horizontal side is not
smaller than the vertical one and is $l$ long. Then there are no more than $l$
squares
$y_{l1},\ldots,y_{ll}$,
in the $l$th column of this rectangle, we shall
apply the elements from $\fn_+^{(1)}$ that shift the $l$th column to the right
by $1,\ldots,(l-1)$ squares (if the vertical side is greater than or equal to
the horizontal one, then the lower row must be shifted upward by the elements
from $\fn_+^{(2)}$). We  ragard as variables only
$y_{l1},\ldots,y_{ll}$, i.e., we assume that the other variables have
{\it general\/} values. We search for a nontrivial {\it linear\/} expression
from $y_{l1},\ldots,y_{ll}$ which are vanish under the corresponding
$(l-1)$ vector fields. Let $\overline\xi=\sum\limits_{i=1\ldots l}a_iy_{i,l}$
be this expression, and then
$$
\left.
{\aligned
\sum\limits_{i=1\ldots l}a_i{\overline y}_{i,l-1}=0\\
-\ -\ -\ -\ -\ -\\
\sum\limits_{i=1\ldots l}a_i,{\overline y}_{i,1}=0
\endaligned}
\right\},
\tag 1
$$
where $a_i$ are unknowns and ${\overline y}_{i,j}$ are $y_{i,j}$ regarded
as constants of the general position. It is clear that then the solution
$\{a_i\}$ is unique, and therefore,
$$
{\overline \xi}=\sum_{i=1\ldots l}a_iy_{i,l}=
\Det\pmatrix
y_{1,l} & {\overline y}_{1,l-1} & \ldots & {\overline y}_{1,1}\\
\vdots & \vdots & & \vdots \\
y_{l,l} & {\overline y}_{l,l-1} & \ldots & {\overline y}_{l,1}
\endpmatrix,
\tag 2
$$
since the minors with respect to the $l$th column taken with the signs
satisfy~(1).

Next, all (and not only linear) expressions from
$\{ y_{1,l},\ldots,y_{l,l}\}$, that vanish under the corresponding
$(l-1)$ vector fields lie in
$$
\C \left[ {\Det}\pmatrix
y_{1,l} & {\overline y}_{1,l-1} & \ldots & {\overline y}_{11}\\
\vdots & \vdots & & \\
y_{l,l} & {\overline y}_{l,l-1} & \ldots & {\overline y}_{l,1}
\endpmatrix \right] = A.
\tag 3
$$
Indeed, let us consider an $l$-dimensional vector space
$V=\{ y_{1l},\ldots ,y_{l,l}\}$.
For every point of $V$ there is a hyperplane going through it and
vanishing under these $l-1$ vector fields.
We draw in $V$ a straight line intersecting these hyperplanes
transversally. Clearly, the function which is vanishing under
the fields is uniquely defined by its restriction to this straight line.
However, (2) defines a linear function of the parameter of this straight line
and any polynomial of the parameter of the straight line lies in $A$ (see
(3)). Then we seek the
$gl_n^{(1)}\oplus gl_n^{(2)}$-singular {\it weight\/} vectors, and any
vector of this kind has in accordance with that has been proved already,
the form $C\cdot\Det_l^k\cdot v$, where $C$ is expression from the least
square. We can apply the preceding arguments to the expression $C$ and
obtain Lemma~1.
\enddemo
\subsubhead
1.4
\endsubsubhead
The next step is to find when $e_0$ acts by zero on the expression
${\Det}^{k_1}_1\cdot \ldots \cdot {\Det}^{k_n}_n\cdot v$
for finding the $gl_{2n}$-singular vectors.
We shall first calculate
$[e_0,{\Det}_k]\ (e_0=x_{11})$.

Let, as in Subsec~1.0, $A_k$ be the matrix
$(y_{ij})\ i,j=1\ldots k$. We denote by
$\widetilde{A_k}$ the matrix
$(y_{ij})\ i,j=2\ldots k$. This is a
$(k-1)\times (k-1)$ matrix. Suppose, furthermore, that
$z^+_i$  is an element from
${\frak n}_-^{(1)}$ lying in the $i$th column above
$y_{1i}$ and $z^-_j$  is an element from
${\frak n}_-^{(2)}$ lying in the $j$th row on the right of
$y_{j1}$. We write
${\Det}_k={\Det}A_k$  first with respect to the upper row and then with
respect to the right column. We obtain
$$
{\Det}A_k=(-1)^{1+k}y_{11}\cdot {\Det}
\widetilde{A_k}+\sum_{i,j=2\ldots k}(-1)^{i+j}y_{1i}y_{j1}\cdot
{\Det}\widetilde{A_{ij}},
$$
where $\widetilde{A_{ij}}$ is the matrix $\widetilde{A_k}$ without the
$i$th row. The matter is that only $y_{1i}$ and $y_{j1}$ do not commute
with $e_0$.

We have
$$
[e_0,{\Det}_k]=(-1)^{1+k}{\Det}\widetilde{A_k}\cdot
\alpha_0^{\vee}+\sum_{i,j=2\ldots k}(-1)^{i+j}{\Det}
\widetilde{A_{ij}}(z^+_i\cdot y_{j1}-y_{1i}z^-_j).
$$
Since $[z^+_i,y_{j1}]=-y_{ji}$, we have
$$
\align
[e_0,{\Det}_k]&=(-1)^{1+k}{\Det}\widetilde{A_k}\cdot
\alpha_0^{\vee}+\sum_{i,j=2\ldots k}(-1)^{i+j-1}\cdot y_{ji}
{\Det}\widetilde{A_{ij}}+\\
&\sum_{i,j=2\ldots k}(-1)^{i+j}\cdot {\Det}
\widetilde{A_{ij}}(y_{j1}z^+_i- y_{1i}\cdot z^-_j).
\tag 4
\endalign
$$
However,
$$
\sum_{i,j=2\ldots k}(-1)^{i+j-1}y_{j1}\cdot
\widetilde{A_{ij}}=(-1)^{k+1}\cdot (k-1)\cdot {\Det}\widetilde{A_k}.
\tag 5
$$
Therefore from (4) and (5) we have
$$
\align
[e_0,{\Det}_k]&=(-1)^{1+k}\cdot {\Det}\widetilde{A_k}\cdot
(\alpha_0^{\vee}+k-1) + \\
&\sum_{i,j=2\ldots k}(-1)^{i+j}\cdot {\Det}\widetilde{A_{ij}}
(y_{j1}\cdot z^+_i- y_{1i}z^-_j).
\tag 6
\endalign
$$

It is clear that the second term in~(6) belongs to ${\Cal U}({\frak
A}_-)\cdot {\frak n}_-^{(1)}\oplus  {\Cal U}({\frak
A}_-)\cdot {\frak n}_-^{(2)}$ and therefore acts by zero on the highest
weight vector $v$ in $Ind_{\mu}$.  We have proved the following lemma.
\proclaim{Lemma~2}
${\Det}_k\cdot v$ is a singular vector in
$Ind_{\mu}$ for $\mu =-k+1$.
\endproclaim

Obviously, ${\Det}_1^{k+1}\cdot v$ is a singular vector in
$Ind_{\mu}$ for $\mu =k\geq 0$.  We have thus proved that for $\mu
\in \Z_{\geq -n+1}$ the representations $Ind_{\mu}$ are reducible.

Note the following important consequence: let $\mu =-k+1+t$, where $t$
is small. It follows from~(6) that $e_0\cdot {\Det}_k\cdot v$
is a {\it first}-order infinitesimal. This is also true in the case of a
singular vector of the form ${\Det}_1^{k+1}\cdot v$ for $\mu =k\geq 0$.
\subsubhead
1.5
\endsubsubhead
Here we calculate $[e_0,{\Det}_k^l]$.
\proclaim{Lemma~3}
$$
\align
[z_i^+, {\Det}_k]=0, & \\
[z_i^-, {\Det}_k]=0, & (i=2\ldots k).
\endalign
$$
\endproclaim
\demo
{Proof}
$ad(z_i^+)$ acts by a vector field on $S^*({\frak A}_-)$
which shifts the first column to the $i$th one and  $ad(z_i^-)$ acts by a
vector field which shifts the first row to the $i$th one.
Hence, for $i\leq k$, we obtain the determinant of a matrix with two similar
columns (rows).
\enddemo

\proclaim
{Lemma~4}
$$
[e_0,{\Det}_k^l]=(-1)^{1+k}\cdot l\cdot \biggl({\Det}\widetilde{A_k}\cdot
{\Det}_k^{l-1}\biggr) (\alpha_0^{\vee}+k-l)+x,
$$
where $x\in {\Cal U}({\frak A}_-)\cdot {\frak n}_-^{(1)}\oplus
{\Cal U}({\frak A}_-)\cdot
{\frak n}_-^{(2)}$.
\endproclaim
\demo{Proof}
The weight of  ${\Det}_k$ is equal to
$-\biggl( k\alpha_0+(k-1)(\alpha_{-1}+\alpha_1)+\ldots
+(\alpha_{-k+1}+\alpha_{k-1})\biggr) =\Theta_k;\
\Theta_k(\alpha_0^{\vee})=-2$. According to Lemma~3 and relation~(6), we have
$$
\align
[e_0,{\Det}_k^l]&=(-1)^{1+k}\cdot \biggl[ {\Det}\widetilde{A_k}\cdot
(\alpha_0^{\vee}+k-1)\cdot {\Det}_k^{l-1}+\\
&{\Det}_k\cdot {\Det}\widetilde{A_k}\cdot
(\alpha_0^{\vee}+k-1)\cdot {\Det}_k^{n-2}+\ldots \biggr] +x,
\endalign
$$
where $x\in {\Cal U}({\frak A}_-)({\frak n}_-^{(1)}\oplus {\frak n}_-^{(2)})$.
$$
(\alpha_0^{\vee}+k-1)\cdot {\Det}_k={\Det}_k\biggl(
\alpha_0^{\vee}+k-1+\Theta_k(\alpha_0^{\vee})\biggr) =
{\Det}_k(\alpha_0^{\vee}+k-3).
$$
Therefore,
$$
\align
[e_0,{\Det}_k^l]&=(-1)^{1+k}\biggl( {\Det}\widetilde{A_k}\cdot
{\Det}_k^{l-1}\biggr)\cdot \Biggl(
(\alpha_0^{\vee}+k-1)+(\alpha_0^{\vee}+k-3)+\\
&\ldots +\biggl(
\alpha_0^{\vee}+k-1-2(l-1)\biggr) \Biggr) +x.
\endalign
$$
The part dependent on $\alpha_0^{\vee}$ is equal to
$(\alpha_0^{\vee}+k-l)$.
\enddemo

\proclaim
{Corollary~1}
\rom {(1)} ${\Det}^l_k\cdot v$ is a singular vector in
$Ind_{\mu}$ only for $\mu = l-k$.

\rom {(2)} If $\mu =l-k+t$, where $t\in \C$ is small, then
$[e_0,{\Det}_k^l]\cdot v$ is the first-order infinitesimal with respect to
$t$.
\endproclaim
\subsubhead
1.6
\endsubsubhead
Let $\xi \cdot v={\Det}_1^{l_1}\cdot \ldots \cdot
{\Det}_k^{l_k}\cdot v$ be a singular vector in $Ind_{\mu}$ with $l_k\ne
0$. We denote $\widetilde{\Det}_k=(-1)^{k+1}\cdot l_k\cdot
{\Det}\widetilde{A_k}$, and then
$$
\aligned
&[e_0,{\Det}_1^{l_1}\cdot \ldots \cdot {\Det}_k^{l_k}]=\\
&\widetilde{\Det}_1\cdot {\Det}_1^{l_k-1}{\Det}_2^{l_2}\cdots
{\Det}_k^{l_k}\biggl( \alpha_0^{\vee}+1-l_1-2(l_2+l_3+\ldots
+l_k)\biggr) +\\
&{\Det}_1^{l_1}\widetilde{\Det}_2{\Det}_2^{l_2-1}
{\Det}_3^{l_3}\cdots {\Det}_k^{l_k}\biggl(
\alpha_0^{\vee}+2-l_2-2(l_3+\ldots +l_k)\biggr) +\\
&\dots +{\Det}_1^{l_1}\cdot \cdots \cdot
{\Det}_{k-1}^{l_{k-1}}\cdot \widetilde{\Det}_k\cdot
{\Det}_k^{l_{k-1}}(\alpha_0^{\vee}+k-l_k)+x,
\endaligned
\tag 7
$$
where $x\in {\Cal U}({\frak A}_-)({\frak n}_-^{(1)}\oplus {\frak
n}_-^{(2)})$.

We shall consider all determinants as elements of the ring
$\C[y_{ij};\,i,j=1\ldots n]$. The last term in this ring cannont be divided
by ${\Det}_k^{l_k}$ but can only be divided ${\Det}_k^{l_k-1}$, and therefore,
if $\xi v$ is a singular vector in $Ind_{\mu}$, then either the last term
in~(7) (with the exception for $x$) acts on $v$ by zero or all the others.
It follows from~(7) that the parts in the parentheses, dependent on
$\alpha_0^{\vee}$,
have the form $(\alpha_0^{\vee}-a_s)$, where $a_s$ increases. Since
$l_k\ne 0$ by hypothesis, the last term acts on $v$ by zero.
Hence follows the part of Theorem~1 stating that $\mu \in
\Z_{\geq -n+1}$.  Let us consider the other terms.
The last term cannont be divided by
${\Det}_{k-1}^{l_{k-1}}$
but can only
be divided  by ${\Det}_{k-1}^{l_{k-1}-1}$ and the other terms can be divided
by  ${\Det}_{k-1}^{l_{k-1}}$. Therefore this term acts
on $v$ by zero, and this is impossible since the sequence $a_s$
increases. This proves that the vectors ${\Det}_k^l\cdot v$
exhaust all singular vectors in $Ind_{\mu}$ and also proves Theorem~1.
\subhead
2. Shapovalov's Form on $Ind_{\mu}$
\endsubhead

\subsubhead
2.0
\endsubsubhead
Let us consider any semisimple Lie algebra $g$ and assume that $M$
is a certain Verma module over $g$. Then the contragradient module $M'$
has the same highest weight as $M$, and therefore, there is a unique
$g$-homomorphism $\alpha :M\to M'$. On every level $V_k$ of the module $M$,
the map $\alpha$ defines the bilinear form $\alpha_k: V_k\otimes
V_k\to \C$. The form $\alpha_k$ is known as {\it Shapovalov's form\/} of the
Verma module $M$.  We have $g={\frak n}_-\oplus {\frak h}\oplus
{\frak n}_+$ and, for $f_i\in {\frak n}_-$, denote $e_i=\omega (f_i)$, where
$\omega$ is the Chevalley involution. It is clear that
$\langle f_{i_1}f_{i_2}\cdots
f_{i_k}v,\ f_{j_1}\cdots f_{j_l}v\rangle =e_{i_k}\cdots
e_{i_1}f_{j_1}\cdots f_{j_l}v$.  Next, the form $\alpha_k$
is {\it symmetric}. The vector $v\in V_k$ lies in a certain proper
submodule in $M$ if and only if $\langle v,w\rangle =0$ for
all $w\in V_k$. If Shapovalov's form is nondegenerate on $V_k$,
then the maximal submodule in $M$ does not intersect $\bigoplus\limits_{i\leq
k}V_i$. Let us consider a quotient module $L$ of the module $M$,
$$
L = M\big/ M_1.
$$

Then $(V_k\cap M_1)\bot V_k$, and therefore, there exists Shapovalov's form
$\alpha_k: L\otimes L\to \C$ with the properties given above. In particular,
$\alpha_k$ is nondegenerate if and only if $L$ is irreducible.

Let us consider Shapovalov's form on $Ind_{\mu}$. For $\mu \notin \Z_{\geq
-n+1}$ it is nondegenerate.
The determinant of Shapovalov's form is not uniquely defined. However, if we
regarded this determinant as a function of $\mu$, then the multiplicity of
zero is uniquely defined.

Our aim is to determine the multiplicity of zero of the determinant of
Shapovalov's form on the levels of $Ind_\mu$.
\subsubhead
2.1
\endsubsubhead
Lemma~1 states that
$w={\Det}_1^{k_1}\cdot \ldots \cdot {\Det}_n^{k_n}\cdot v$
is a $gl_n^{(1)}\oplus gl_n^{(2)}$-singular vector in $Ind_{\mu}$ for any
$k_1,\ldots ,k_n\geq 0$ and these monomials exhaust
$gl_n^{(1)}\oplus gl_n^{(2)}$-singular vectors.

\proclaim
{Lemma~5}
For different $(k_1,\ldots ,k_n)$ the $gl_n^{(1)}\oplus gl_n^{(2)}$-weights
of the corresponding monomials $w$ are different.
\endproclaim

\demo{Proof}
The proof is obvious.
\enddemo

This very simple result has a fundamental significance.
\proclaim
{Lemma~6}
$Ind_{\mu}=\oplus$ \rom(of irreducible finite-dimensional
$gl_n^{(1)}\oplus gl_n^{(2)}$-modules with different highest weights\rom{).}
\endproclaim

\demo{Proof}
We have only to prove  that ${\frak n}_-^{(1)}\oplus {\frak
n}_-^{(2)}$ acts on $Ind_{\mu}$ locally finitely. The operators $ad{\frak
n}_-^{(1)} (ad{\frak n}_-^{(2)})$ shift any element $y_{ij}$ to the
left (downward) and the elements of the  $n$th column
(of the $n$th row) commute with ${\frak
n}_-^{(1)} ({\frak n}_-^{(2)})$.
\enddemo

\example{Example} The $gl_n^{(1)}\oplus gl_n^{(2)}$-modules corresponding
to ${\Det}_1\cdot v,\ {\Det}_2\cdot v,\ldots$, ${\Det}_n\cdot v$
are, respectively,
$$
\Lambda^{n-1}V\otimes \Lambda^{n-1}V, \Lambda^{n-2}V\otimes
\Lambda^{n-2}V,\ldots ,\Lambda^1V\otimes \Lambda^1V,\C,
$$
where $V$ is the tautological $n$-dimensional representation of $gl_n$.
\endexample
\subsubhead
2.2
\endsubsubhead
Since the weight of all arising $gl_n^{(1)}\oplus
gl_n^{(2)}$-modules are different, under the mapping $\alpha : Ind_{\mu}\to
Ind'_{\mu}$ that defines Shapovalov's form every irreducible
$gl_n^{(1)}\oplus gl_n^{(2)}$-module $L$ passes into $L'$.
It follows that all irreducible
$\g_n^{(1)}\oplus\g_n^{(2)}$-modules arising in $\I$ are pairwise
orthogonal in the sence of Shapovalov's form.
Thus, if we choose a base on a certain level of $\I$ that agrees with
the decomposition of $\I$ into the direct sum of
$\g_n^{(1)}\oplus \g_n^{(2)}$-modules, then the matrix of Shapovalov's form
is block-diagonal in this base. Therefore, the determinant of Shapovalov's
form is equal to the product of the determinants of these blocks.

Next, let $w=\D_1^{k_1}\cdot \ldots \cdot \D_n^{k_n}\cdot
v$ be a singular vector $\g_n^{(1)}\oplus \g_n^{(2)}$-module $L$.
A typical vector of the intersection of $L$ with a certain level of $\I$
has the form $f_{i_s}\ldots f_{i_1}w$, where $f_{i_1},\ldots ,f_{i_s}\in
\n^{(1)}_-\oplus \n_-^{(2)}$.  The corresponding block in Shapovalov's
matrix consists of elements
$w'e_{j_1}\ldots e_{j_t}f_{i_s}\ldots f_{i_1}w$, where
$e\ldots$ and $w'$ are images of $f\ldots$ and $w$ under the Chevalley
involution. Since
$\n_+^{(1)}$ É $\n_+^{(2)}$ acts trivially on $w$, we find that the
block of Shapovalov's matrix on the $j$th level of $\I$, corresponding to the
$\g_n^{(1)}\oplus \g_n^{(2)}$-module $L$, is equal to
$\langle w,w\rangle^{p_j(L)}\cdot A$, where
$p_j(L)$ is the dimension of the intersection of the $j$th level of $\I$
with $L$ and
$A$ is the matrix of Shapovalov's form of the
$\g_n^{(1)}\oplus \g_n^{(2)}$-module $L$
on this level.  Note that here only
$\langle w,w\rangle$ depends on $\mu$ and
$\det A\ne 0$ since $L$ is irreducible. Thus the problem of determining
the multiplicity of the zero of the determinant of Shapovalov's form
(on any level) reduces to the search for
$p_w(\mu )=\langle w,w\rangle$ for any monomial
$w=\D_1^{k_1}\cdot \ldots \cdot \D_n^{k_n}\cdot v$.
\subsubhead
2.3
\endsubsubhead
The direct calculation of $\langle w,w\rangle$, even for $w=\D_k\cdot
v$ is impossible, and therefore, we do the following: first, we prove
that $\e \ p_w(\mu )=\sum\limits^n_{i=1}i\cdot k_i$ (which was to be
expected), and second, we find possible roots of $p_w(\mu )$ and their
{\it maximal\/} possible multiplicities. It turns out that the sum of
these multiplicities is exactly equal to $\sum\limits^n_{i=1}i\cdot k_i$.
\proclaim{Lemma~7}
$p_w(a)=0\Leftrightarrow w$ belongs to the submodule in $\In$.
\endproclaim

\demo{Proof}
If $w$ belongs to the submodule, then it is clear that
$\langle w,w\rangle =0$.  On the other hand, in the block base (see
Subsec.~2.2) of a certain level of $\I \quad w$ is a $1\times 1$ block
since $w$
is the only element of some $\g_n^{(1)}\oplus \g_n^{(2)}$-module
$L$ on its level.  Therefore, $\langle w,w\rangle =0\Rightarrow \langle
w_1,w\rangle =0$ for any $w_1$ from the corresponding level.
Consequently, $w$ belongs to the proper submodule in $\I$.
\enddemo

It follows from Lemma~7 and Theorem~1 that $p_w(\mu )$ can have only integer
$\geq -n+1$ roots. Let us first consider $w=\D_k\cdot v$.

\proclaim{Lemma~8}
For $w=\D_k\cdot v\quad p_w(\mu )$ can only have the roots
$$
-k+1,\ldots ,-1, 0
$$
\endproclaim
\noindent
(a total of $k$ roots, and $p_w$ is a of degree $k$ polynomial; however,
there can be multiplicities).

\demo{Proof}
It follows from Lemma~7 that if $p_w(a)=0$, then we can obtain a certain
singular  vector in $\In$ by acting on $w$ by operators from
$\n_+\subset \g_{2n}$.

(1) $a\geq -k+1$.

The action of $ad \A_+$ cannont sent the element from the
$k\times k$ square of $A_k$
into $\A_-$, and the first singular vector in $Ind_{-k}$ is
$\D_{k+1}\cdot v$.

(2) $a\leq 0$.

It is true that for every $w=\D_1^{k_1}\ldots \D_n^{k_n}\cdot v$, the action
$\A_+$ cannont increase the degree of $w$ with respect to any $y_{ij}\in \A_-$.
Let us prove this for $w=\D_k\cdot v$ and $y_{11}=f_0$.
\vskip+5cm
\centerline{Fig. 2}

$\D_k$ is the sum of monomials each of which contains exactly one element
from the first row and exactly one element from the first column
(they can coincide). Let $z_+^1,
z_+^2, z_+^3,\ldots$ be elements lying in the first column of $\A_-$
outside  of $\A_-$ and $z_-^1, z_-^2, z_-^3,\ldots$ be elements lying
in the first row of $\A_-$ outside of $\A_-$ (see Fig.~2).
The reader can easily verify that the element
$f_0$ can be obtained in only two ways:

(1) we apply  $ad (x_{is})$ to $y_{1i}$,
obtain $z_-^{s-1}$, and apply $ad(z_-^{s-1})$ to $y_{s1}$;

(2) we apply $ad (x_{si})$ to $y_{i1}$,
obtain $z_+^{s-1}$, and apply $ad(z_+^{s-1})$ to $y_{1s}$.

Thus, one element from the first row and one element from the first
column must be used in the formation of $f_0$.
In general, an element from the $i$th row and an element from the
$j$th column must be used in the  formation
of the element $y_{ij}$. Therefore, in the process of apploication of
elements from $\A_+$ to $\D_k$, every monomial will contain not more than one
element from every row and every column.

For $\mu =1$, the singular vectors $f_0^2,\ldots$ in
$Ind_1$ contain $f_0$ of the $\geq 2$ degree, and therefore,  is
impossible to obtain them from $\D_k$ by the application of elements from
$\A_+$.
\enddemo

The proof of Lemma~8 implies the following corollary.
\proclaim{Corollary~2}

\rom {(1)} For any monomial $w=\D_1^{k_1}\cdot \ldots \cdot
\D_n^{k_n}\cdot v$ the application of elements from $\A_+$
cannont increase the
degree of $w$ with respect to any $y_{ij}\in \A_-$,

\rom {(2)} the roots of $p_w(\mu )$ lie in the interval $[-l+1, (\sum k_i)-1]$,
where $l$ is the largest $i$ for which $k_i\ne 0$.
\endproclaim
\subsubhead
2.4
\endsubsubhead
Here we find an upper estimate for the multiplicities of the roots of
$p_w(\mu )$ given by Corollary~2~(2).

\proclaim{Lemma~9}
If $\mu =a$ is a root of multiplicity $d$ of the polynomial
$p_w(\mu )$, then, for $\mu =a+t$ ($t$ is small), $\langle w_1,w\rangle_{a+t}$
is either equal to 0 or is an infinitesimal of the order of exactly $d$
with respect to $t$. In other words, for any
$e_{j_1},\ldots , e_{j_s}\in \n_+\subset \g_{2n}$,
$(e_{j_s}\ldots e_{j_1}\cdot w)\cdot v$ is either 0 or an infinitesimal of
order $d$.
\endproclaim

\demo{Proof} This is similar to the proof of Lemma~7.
\enddemo

It follows from Lemma~9 that the required multiplicity is equal to the
multiplicity of zero for $\mu=a+t$ in the expression
$(e_{j_s}\ldots e_{j_1}\cdot w)\cdot v$ for
{\it all\/} $e_{j_1},\ldots , e_{j_s}\in \n_+\subset \g_{2n}$.

\proclaim{Lemma~10}  The maximal possible multiplicity of the root $\mu =a$
in $p_w(\mu )$ is not larger than the number of singular vectors in
$\In$ lying at the levels not exceeding the level of $w$.
\endproclaim

\demo{Proof}
According to Corollary~1~(2), every singular vector increases the order
with respect to $t$ by unity.
On the other hand, if ${\widetilde w}\in \In$
is not a singular vector, then $e_j {\widetilde w}$ in $Ind_{a+t}$
is not an infinitesimal and not zero for some $e_j\in \n_+\subset \g_{2n}$.
Iterating this process we derive the statement of the lemma from Lemma~9.
\enddemo

In conjuction with Corollary~2 (1) we get

\proclaim{Lemma~11} The maximal posiible multiplicity of the root
$\mu =a$
in $p_w(\mu )$ is equal to the number of singular vectors in $\In$ lying
not lower than $w$ and having \rom(in the notation $\D^r_p\cdot v$\rom)
a degree with respect to any
$y_{ij}\in \A_-$ not higher than $w$.
\endproclaim

\remark{Examples} In Subsec.~2.5 we shall prove that $\e \ p_w(\mu)=
\sum\limits_{i=1\ldots n}i\cdot k_i$.

(i) $w=\D_k\cdot v$.

For $-k+1\leq a\leq 0$, the conditions of Lemma~11 are satisfied only
by the singular vector $\D_{-a+1}\cdot v$, and therefore, every one
of the possible $k$ roots has a multiplicity not higher than~1, and
since $\e \ p_w(\mu )=k$, we have
$p_w(\mu )=\mu \cdot (\mu +1)\cdot \ldots \cdot
(\mu +k-1)$ with an accuracy up to a constant.

(ii) $w=\D_2\cdot \D_3\cdot v$.

The possible roots of $p_w(\mu )$ are
$$
+1,\ 0,\ -1,\ -2
$$

for $\mu =-2$, singular vectors in $\I$ are $\D_3\cdot v, \D_4^2\cdot
v,\ldots$; the conditions of Lemma~11 are satisfied only by $\D_3\cdot v$,
the maximal multiplicity is~1;

for $\mu =-1$ singular vectors are $\D_2\cdot v, \D_3^2\cdot v,\ldots$;
the conditions of Lemma~11 are satisfied only by $\D_2\cdot v$; the maximal
multiplicity is~1;

for $\mu =0$, singular vectors are $\D_1\cdot v, \D_2^2\cdot v, \D_3^3\cdot
v,\ldots$;
the conditions of Lemma~11 are satisfied only by
$\D_1\cdot v, \D_2^2\cdot v$;
the maximal multiplicity is 2;

for $\mu =+1$, singular vectors are $\D_1^2\cdot v, \D_2^3\cdot v,\ldots$;
the conditions of Lemma~11 are satisfied only by
$\D_1^2\cdot v$;  the maximal multiplicity is 1.

On the other hand, $\e \ p_w=5$, and therefore
$$
p_w(\mu ) = (\mu -1)\cdot \mu^2\cdot (\mu +1)(\mu +2).
$$
with the accuracy up to a constant.
\endremark
\proclaim{Lemma~12} For $w=\D_1^{k_1}\cdot \ldots \cdot
\D_n^{k_n}\cdot v$ we have
$$
\e \ p_w(\mu ) = \sum_{i=1\ldots n}i\cdot k_i.
$$
\endproclaim

The proof will be given in Subsec~2.5.

The situation encountered in the examples given above is the same for any $w$.
Indeed, we can use another technique to calculate the sum of the maximal
possible multiplicities. We can obtain
$\D_n^{k_n}$, $\D_n^{k_n-1},\ldots$,\linebreak
$\D_n$; $\D^{k_n}_{n-1},\ldots ,
\D_{n-1};\ldots$; $\D_1^{k_n},\ldots , \D_1$,
from $\D_n^{k_n}$, certainly, for different
$\mu$. We have the total of $n\cdot k_n$ singular vectors.
Repeating this reasoning for $\D_i^{k_i}$ for any $i$, we find
that the sum of the maximal possible multiplicities is
$\sum\limits_{i=1\ldots n}i\cdot
k_i$. We have proved

\proclaim{Lemma~13} The sum of multiplicities calculated with the use
of Lemma~11 is
$\sum\limits_{i=1\ldots n}i\cdot k_i$.
\endproclaim

It is clear that we can find $p_w(\mu )$ (with an accuracy to
within the multiplication by a constant) for any
$w=\D_1^{k_1}\cdot \ldots \cdot
\D_n^{k_n}\cdot v$.

There are two ways of formulating the answer in a concise form (the second
form is the main one), namely,

\remark{Technique~1}
For every $\D_k$, we must write the corresponding polynomial
$p_k(\mu )=\mu (\mu +1)\ldots (\mu +k-1)$ and shift this $\D_k$
to the right, to the end of the monomial; when shifting by one, the
determinant (of the first degree) $p_k(\mu )$
is replaced by $p_k(\mu -1)$. We must shift all the polynomials to the end
of the monomial and then multiply them. (We must do this in the order
$$
w = \D_1^{k_1}\cdot \ldots \cdot \D_n^{k_n}\cdot v\text{)}.
$$
\endremark

\remark{Technique~2}
If we consider $w$ to be the highest weight of the corresponding
$\g_n^{(2)}$-module, then a certain Young diagram ${\Cal D}(w)$ corresponds
to it (here $n\gg \e \ w$). We must cut ${\Cal D}(w)$
by oblique straight lines under the angle of $45^{\circ}$ as is shown in
Fig.~3.
\topinsert
\vskip+5cm
$$
\gather
w=\D_2\cdot \D_3\cdot v\\
p_w(\mu )=(\mu -1)\mu^2(\mu +1)(\mu +2)
\endgather
$$
\centerline{Fig.~3}
\endinsert
\endremark

The number of squares on a certain diagonal is the multiplicity
of the corresponding root, from the smaller to the larger roots if we
cut upward. We denote by $l({\Cal D})$ the number of diagonals and by
$n_1,\ldots ,n_{l(\Cal D)}$ the number of squares on the diagonals
if we cut upward. Then
$$
p_w(\mu ) = \prod^{l(\Cal D)}_{i=1}(\mu +k({\Cal D})-i)^{n_i},
\tag 8
$$
where $k({\Cal D})$ is the maximal height of ${\Cal
D}(w)$ ($p_w(\mu )$ with an accuracy to within the multiplication
by a constant). Thus,
the length of the ``central'' diagonal that passes through the upper
left vertex of ${\Cal D}(w)$ is rqual to the multiplicity of the root
$\mu =0$.
\subsubhead
2.5
\endsubsubhead
Here we prove Lemma~12. As we have mentioned, the direct calculation of
$p_w(\mu )$ is impossible, even for $w=\D_k\cdot v$.
However, it is very easy to calculate explicitly the highest coefficient and
prove that it is not equal to zero.

Let $w=\D_1^{k_1}\cdot \ldots \cdot \D_n^{k_n}\cdot v$. We shall calculate
$\langle w,w\rangle$. We denote by $p^i_l$ the monomials, the terms of
$\D_l$, taken with the signs
$(i=1,\ldots ,l!)$ and by $q^i_l$ the images of $p^i_l$ under the Chevalley
involution.  The sign that appears due to the Chevalley involution does not
depend on the choice of a definite monomial but only depends on $l$. Only the
terms
$$
\biggl( (q_n^{i_{k_n}}\ldots q_n^{i_1})(q_{n-1}^{j_{k_{n-1}}}\ldots
q_{n-1}^{j_1})\ldots (p_{n-1}^{j_{k_{n-1}}}\ldots p_{n-1}^{j_1})
(p_n^{i_{k_n}}\ldots p_n^{i_1})\biggr) \cdot v,
\tag 9
$$
contribute to the highest monomial namely, the terms which are ``symmetric''
with respect to $p$ and $q$. Due to this symmetry, the signs of all monomials
are cancelled out. Now we set
${\widetilde q}_l^{\ \sigma}=\prod\limits_{i=1\ldots l}x_{\sigma (i)i}$ for
$p^{\sigma}_l=\prod\limits_{i=1\ldots
l}y_{i\sigma (i)}$. Then, when we replace
$q$ by $\widetilde q$, relation~(9) remains valid with an accuracy to within
the common sign dependent on $w$. Now we shift $\widetilde q$ to the right and,
as a result, no signs appear. We have proved the lemma.
\subhead
3. The Jantzen Filtration: Formulas for the Characters
of Irreducible Representations and Combinatorial Identities
\endsubhead
\subsubhead
3.0
\endsubsubhead
The considerations of Sec.~2 result in numerous formulas and identities.
Indeed, let $a\in \Z_{\geq -n+1}$ be a critical value of the parameter,
and let $\In^{(1)}$ be the kernel of Shapovalov's form on
$\In$. Then $\In^{(1)}$ is a submodule and $\In \big/ \In^{(1)}$ is
irreducible. Thus, from the results of Sec.~2 we can deduce formulas
for the
characters of irreducible representations of $\g_{2n}$ and ${\widehat
\g}_{\infty}$. For example, at $\mu =1$ for
${\widehat \g}_{\infty}$ the character of the irreducible
representation of ${\widehat \g}_{\infty}$
with zero highest weight and the central charge~1 is equal, by the Weyl
formula, to
$\prod\limits_{i\geq 1}{\displaystyle\frac{1}{(1-q^i)}}$,
and from the results of Sec.~2 the same character is equal (see
Subsec.~3.2) to
$1+\sum\limits_{k\geq1}
{\displaystyle\frac{q^{k^2}}{(1-q)^2\ldots (1-q^k)^2}}$.

Thus we obtain the Euler identity
$$
\prod_{i\geq 1}\frac{1}{(1-q^i)} = 1 + \sum_{k\geq
1}\frac{q^{k^2}}{(1-q)^2\ldots (1-q^k)^2}.
\tag 10
$$

If we set $\mu =2,3,\ldots$, then we get the ``higher analogs''
of identity~(10). On the other hand, for $\mu =-1,-2,\ldots$, we
get relations for the character of the corresponding irreducible
representation $\I\big/ \I^{(1)}$ (see (17), (18) belolw).
Note that by considering
$\g_{2n}$ rather than ${\widehat \g}_{\infty}$, we obtain the
``finite form''
of identity~(10) and of the highest identities (see Subsec.~3.2).
\subsubhead
3.1. Definition
\endsubsubhead
For $k\in \Z_{\geq 1}$, we set $\In^{(k)}=\{
w\in \In \ \vert \ \langle w, w_1\rangle_{a+t}$ is divisible by $t^k$
for any $w_1\in \In \}$. $\In^{(k)}$ is the {\it $k$th term of the Jantzen
filtration}. $\In^{(k)}$
is a submodule in $\In$ for any integer $k\geq 1$, and $\In \supset
\In^{(1)}\supset \In^{(2)}\supset \ldots$ is the
{\it Jantzen filtration\/} of the module $\In$.

\proclaim{Lemma~14
\rm(A description of the Jantzen filtration of $\In$)} Let
$a$ be a critical value of the parameter. Then  $\In^{(k)}$ is the union
of all irreducible $\g_n^{(1)}\oplus \g_n^{(2)}$-modules, the scalar
square $p_w(\mu )$ of whose highest weight vector $w$ has
$\mu =a$ a root of multiplicity $\geq k$.
\endproclaim
\demo{Proof}
According to Lemma~9, the condition imposed on $w$
formulated in the lemma is equivalent to the fact that $w\in \In^{(k)}$.
However, since $\In^{(k)}$ is a
$\g_{2n}$-submodule, $\In^{(k)}$ contains the whole irreducible
$\g_n^{(1)}\oplus \g_n^{(2)}$-module with the weight vector
$w$.
\enddemo

We have proved the following relationship (that always exists) of the
Jantzen filtration with zero multiplicity of the determinant of
Shapovalov's form.

\proclaim{Corollary~3}
Let $V_j$ be $j$th level of
$\In$. Then the multiplicity of zero \rom(with respect to $\mu$ at the point
$\mu =a$\rom) of the determinant of Shapovalov's form on the  $j$th level is
$$
\dim (V_j\cap \In^{(1)}) + \dim (V_j\cap \In^{(2)})+\ldots
$$
\rom(the sum is finite\rom).
\endproclaim

\proclaim{Theorem~2 \rm(Jantzen's conjecture)}
The consecutive quotient modules of the Jantzen filtration are semisimple.
\endproclaim

\demo{Proof}
Due to the existence of the $\g_n^{(1)}\oplus \g_n^{(2)}$-action on
$\In$, as distinct from the case of the Verma modules, in our case, Jantzen's
conjecture can be proved by means of elementary reasoning.

The restriction of Shapovalov's form to $\In^{(k)}$ is divisible by
$t^k$. We divide it by
$t^k$ and then set $t=0$. Then this form induces an invariant symmetric form on
$\In^{(k)}\big/ \In^{(k+1)}$ with values in
$\C$.  We denote it by $\langle ,\rangle_k$. {\it The form $\langle
,\rangle_k$ on $\In^{(k)}\big/ \In^{(k+1)}$ is nondegenerate}. Indeed,
assume that we consider the $l$th level $V_l$ of the module $\In$, and let
$\xi \in V_l$ belong to the kernel of $\langle ,\rangle_k$.
Then $\xi =\xi_1+\ldots
+\xi_s$, where $\xi_1,\ldots ,\xi_s$ are components of $\xi$ relative
to the decomposition  of $V_l$ into the direct sum of the intersections
of $V_l$ with irreducible $\g_n^{(1)}\oplus \g_n^{(2)}$-module
$V_l^{(1)},\ldots ,V_l^{(s)}$.

Assume, for example, that  $\xi_1\ne 0$. It follows from~2.2 that
$V_l^{(1)}\perp \xi_2,\ldots , V_l^{(1)}\perp \xi_s$. Let $w$ be the
highest weight vector of the $\g_n^{(1)}\oplus \g_n^{(2)}$-module
corresponding to $V_l^{(1)}$. Then, in accordance with Lemma~14,
$\langle w,w\rangle$ is vivisible by $t^k$ and not divisible by
$t^{k+1}$, and therefore, the form $\langle ,\rangle_k$ on this
$\g_n^{(1)}\oplus \g_n^{(2)}$-module is nondegenerate. Therefore, $\langle
\xi_1,\xi'_1\rangle_k\ne 0$ for a certain $\xi'_1\in V_l^{(1)}$, and
consequently $\langle \xi ,\xi'_1\rangle_k\ne 0$.

Suppose now that $N\subset \In^{(k)}\big/ \In^{(k+1)}$ is a certain
$\g_{2n}$-submodule. To prove the theorem, it is sufficient to show that
$N\cap N^{\perp}=0$ ($N^{\perp}$ in the sense of $\langle ,\rangle_k$). Let
$\xi \in N\cap N^{\perp}$ and $\xi \ne 0$. Assume, as above, that $\xi_1\in
V_l^{(1)}$ is nonzero. Being a $\g_n^{(1)}\oplus \g_n^{(2)}$-module, $N$
contains all $V_l^{(1)}$, and the reasoning that proves the nondegeneracy of
$\langle ,\rangle_k$ allows us to find $\xi'_1\in V_l^{(1)}\subset N$
such that $\langle \xi_1,\xi'_1\rangle =\langle \xi ,\xi'_1\rangle$
is nonzero, and therefore, $\xi$ is not orthogonal to $N$.
\enddemo
\proclaim{Theorem~3} The submodules generated by singular vector in
$\In$ constitute the complete list of submodules in $\In$. They are embedded
into each other and the
$k$th submodule coincides with $\In^{(k)}$. The consecutive quotient modules
$\In^{(k)}\big/ \In^{(k+1)}$ are simple.
\endproclaim
\demo{Proof}

(1) {\it There is exactly one singular vector in\/}
$\In^{(k)}\big/ \In^{(k+1)}$. In fact only some highest weight vector
of the irreducible $\g_n^{(1)}\oplus
\g_n^{(2)}$-module can be such a singular vector. Now the assertion follows
from Corollary~2 (1) of Sec.~2.

(2) It follows from~(1) and Theorem~2 that the consecutive quotient modules of
the Jantzen filtration are simple. Indeed, assume that
$\D_p^qv$ is a singular vector in
$\In$ and that it belongs to $\In^{(k)}$. Then $\D^{q+1}_{p+1}v$ is the
next singular vector and it belongs to $\In^{(k+1)}$. $\D^q_{p+1}v$
belongs to $\In^{(k)}$ and, since the representation
$\In^{(k)}\big/ \In^{(k+1)}$ is simple,
$$
\D^q_{p+1}\cdot v={\Cal U}(\A_+)\cdot \D_p^q\cdot v+\alpha,
$$
where $\alpha \in \In^{(k+1)}$. It follows from Lemma~14 and Lemma~11
that $\alpha=0$ from the grading considerations. Therefore,
$\D^q_{p+1}\cdot v$ belongs to the submodule in $\In$ generated by the
singular vector $\D^q_p\cdot v$, and this implies that the singular vector
$\D^{q+1}_{p+1}\cdot v$ also belongs to this module. This proves Theorem~3.
\enddemo
\subsubhead
3.2
\endsubsubhead
Here we consider formulas following from the fact that
$\In \big/ \In^{(1)}$ is an irreducible module for $a>0$.

First of all, let us consider a representation of ${\widehat \g}_{\infty}$ (or
of $\g_{\infty ,\text{ fin}}$) with a zero highest
weight and the central charge
$\mu =1$ (or with one zero label respectively). Then $\D_1^2\cdot
v$ is the first singular vector, and it follows from Lemma~11 that outside of
$Ind_1^{(1)}$ there are $\g^{(1)}_{\frac{\infty}{2}}\oplus
\g^{(2)}_{\frac{\infty}{2}}$-modules with the
highest weight vectors $v, \D_1\cdot v, \D_2\cdot v,\ldots$ and
only these vectors.
It is clear that the character of the irreducible
$\g^{(2)}_{\frac{\infty}{2}}$-module with the highest vector
$\D_k\cdot v$ is equal to ${\displaystyle\frac{1}{(1-q)\ldots
(1-q^k)}}$, the character of the
$\g^{(1)}_{\frac{\infty}{2}}\oplus \g^{(2)}_{\frac{\infty}{2}}$-module is
equal to its square ${\displaystyle\frac{1}{(1-q)^2\ldots (1-q^k)^2}}$, and
the weight of $\D_k$ is equal to $k^2$. Hence,
$$
\frac{1}{\prod\limits_{i\geq 1}(1-q^i)} = 1 + \sum_{k\geq
1}\frac{q^{k^2}}{(1-q)^2\cdot \ldots \cdot (1-q^k)^2}.
$$
(relation~(10), the Euler identity).

Let us write out the``finite form'' of the last identity, replacing
${\widehat\g}_{\infty}$ by $\g_{2n}$: only $v, \D_1\cdot v,\ldots , \D_n\cdot
v$ remain. The character of the irreducible $\g_n^{(2)}$--module
with the highest weight $\D_k$ is equal to
$$
\frac{(1-q^{n-k+1})\cdot \ldots \cdot (1-q^n)}{(1-q)\ldots (1-q^k)},
\tag 11
$$
and the character of the corresponding irreducible $\g_{2n}$-module is equal to
$$
\frac{(1-q^{n+1})\cdot \ldots \cdot (1-q^{2n})}{(1-q)\ldots (1-q^n)}.
\tag 12
$$
We have
$$
\gathered
{\displaystyle\frac{(1-q^{n+1})\cdot \ldots \cdot (1-q^{2n})}{(1-q)\ldots
(1-q^n)}} =\\
1 +
\sum\limits_{k=1\ldots n}{\displaystyle\frac{q^{k^2}\cdot
(1-q^{n-k+1})^2\cdot \ldots \cdot (1-q^n)^2}{(1-q)^2\cdot \ldots \cdot
(1-q^k)^2}}\\
(n\in \Z_{>0}).
\endgathered
\tag 13
$$

In what follows, we shall not write out the finite forms of the formulas and
identities.

For $\mu =l\geq 1$, the first singular vector in the ${\widehat
\g}_{\infty}$-module  $\I$ is $\D_1^{l+1}\cdot v$, and therefore, the highest weight
vectors of the $\g^{(1)}_{\frac{\infty}{2}}\oplus
\g^{(2)}_{\frac{\infty}{2}}$-modules which lie outside of
$\I^{(1)}$ are the monomials $w=\D_1^{k_1}\cdot \ldots \cdot \D_s^{k_s}\cdot
v$ with $k_1+\ldots +k_s\leq l$. Let ${\Cal D}(w)$ be the corresponding
Young diagram, and let $\chi (w)$ be the character of the corresponding
irreducible
$\g^{(2)}_{\frac{\infty}{2}}$-module. Then
$$
\gather
{\displaystyle\frac{1}{(1-q)(1-q^2)\ldots (1-q^l)^l(1-q^{l+1})^l\cdot
\ldots}}=\\
1 + \sum\limits_{\{ w\vert k_1+\ldots +k_s\leq l\}}q^{k_1+4k_2+\ldots
+s^2k_s}\cdot (\chi (w))^2.
\tag 14
\endgather
$$

We write this explicitly for $l=2$. Here either $w=\D_i\cdot v$ or
$w=\D_i\cdot \D_j\cdot v$ for $i\leq j$. In the second case, the diagram
${\Cal D}(w)$ is shown in Fig.~4:
\vskip+5cm
\centerline{Fig.~4}
$$
\chi (\D_i\cdot \D_j) = \frac{(1-q^{j-i+1})\cdot \ldots \cdot
(1-q^j)}{\prod\limits_{k=1\ldots i}(1-q^k)\cdot \prod\limits_{k=1\ldots
j+1}(1-q^k)}.
$$

We have
$$
\gather
{\displaystyle\frac{(1-q)}{\prod\limits_{s\geq 1}(1-q^s)^2}} = 1 +
\sum\limits_{k\geq 1}{\displaystyle\frac{q^{k^2}}{(1-q)^2\cdot \ldots
\cdot (1-q^k)^2}} +\\
 \sum\limits_{j\geq i\geq 1}{\displaystyle\frac{q^{i^2+j^2}\cdot
(1-q^{j-i+1})^2\cdot \ldots \cdot (1-q^j)^2}{\bigl[
(1-q)^2(1-q^2)^2\cdot \ldots \cdot (1-q^i)^2\bigr] \cdot \bigl[
(1-q)^2\cdot \ldots \cdot (1-q^{j+1})^2\bigr]}}.
\tag 15
\endgather
$$

Finally,
$$
\frac{1}{\prod\limits_{i\geq 1}(1-q^i)^i} = 1 + \sum_{\text{over all}\
w}q^{k_1+4k_2+\ldots +s^2k_s}\cdot \bigl( \chi (w)\bigr)^2,
\tag 16
$$
where $w=\D_1^{k_1}\cdot \ldots \cdot \D_s^{k_s}\cdot v$ and $\chi (w)$
is the character of the irreducible $\g^{(2)}_{\frac{\infty}{2}}$-module
with the corresponding Young diagram. Note that in~(16) the summation is
carried out over {\it all\/} Young diagrams.
\subsubhead
3.3
\endsubsubhead
Let us consider the case $\mu =-1, -2, -3,\ldots$.  For example, if $\mu
=-1$, then $w$ which lie in $Ind_{-1}\setminus Ind_{-1}^{(1)}$ are
$\D_1\cdot v, \D_1^2\cdot v, \D_1^3\cdot v,\ldots$. We denote by
$\chi_{-1}$ the character of the irreducible
${\widehat \g}_{\infty}$-module with the zero highest weight and the central
charge $-1$. We have
$$
\chi_{-1} = 1 + \sum_{k\geq 1}q^k\cdot \frac{1}{(1-q)^2\cdot \ldots \cdot
(1-q^k)^2}.
\tag 17
$$
For $\mu =-l$, the first singular vector in $Ind_{-l}$ is $\D_{l+1}$,
and therefore, the general case is as follows. We denote
$W_l=\{ w=\D_1^{k_1}\cdot
\D_2^{k_2}\cdot \ldots \cdot \D_l^{k_l};\ k_i\geq 0\}$.
Then
$$
\chi_{-l} = \sum_{w\in W_l}q^{k_1+4k_2+\ldots +l^2\cdot k_l}\cdot \bigl(
\chi (w)\bigr)^2.
\tag 18
$$
\subsubhead
3.4
\endsubsubhead
Here we use the whole Jantzen filtration, not only its first term.
According to Theorem~3, the module $\In^{(k)}\big/ \In^{(k+1)}$
is irreducible.
\subsubhead
3.4.1
\endsubsubhead
Let $-k\in \Z_{\leq 0}$, and then the singular vectors in
$Ind_{-k}$ are $\D_{k+1}\cdot v$, $\D^2_{k+2}\cdot v$, $\ldots$.
It is clear that
$Ind_{-k}^{(l)}\big/ Ind_{-k}^{(l+1)}$ is an irreducible module
with the highest weight vector $\D^l_{k+l}\cdot v\ (l\geq 1)$.
We denote the corresponding weight by
$\Lambda^-_l$. We have
$$
\left.
\aligned
&\Lambda^-_l(\alpha^{\vee}_0) = -k - 2l  \\
&\Lambda^-_l(\alpha^{\vee}_{\pm (k+l)}) =l, \qquad l\geq 1\\
&\text{the other}\ \Lambda^-_l(\alpha_i^{\vee}) = 0.
\endaligned
\right\}
\tag 19
$$
The highest weight vectors of the
$\g^{(1)}_{\frac{\infty}{2}}\oplus\g^{(2)}_{\frac{\infty}{2}}$-modules
$w$ lying in $Ind_{-k}^{(l)}\setminus Ind_{-k}^{(l+1)}$ are
$$
W^-_l = \biggl\{ \D_1^{k_1}\cdot \D_2^{k_2}\cdot \ldots \cdot v\
\vert \ \sum_{s\geq k+l}k_s\geq l,\quad \text{and}\ \sum_{s\geq
k+l+1}k_s\leq l\biggr\} .
$$
Then, as usual,
$$
\chi_{\infty}(\Lambda^-_l) = \sum_{w\in
W^-_l}q^{-l(k+l)^2+k_1+4k_2+\ldots}\cdot \bigl( \chi (w)\bigr)^2.
\tag 20
$$
\subsubhead
3.4.2
\endsubsubhead
Let $k\in \Z_{\geq 0}$. The singular vectors in $Ind_k$
are $\D_1^{k+1}\cdot v, \D_2^{k+2}\cdot v,\ldots$. The highest weight
vector of $Ind_k^{(l)}\big/ Ind_k^{(l+1)}$ is $\D_l^{k+l}\cdot v$.
The corresponding highest weight
$$
\left.
\aligned
&\Lambda^+_l(\alpha^{\vee}_0) = -k - 2l  \\
&\Lambda^+_l(\alpha^{\vee}_{\pm l}) =k+l, \qquad l\geq 1\\
&\text{the other}\ \Lambda^+_l(\alpha_i^{\vee}) = 0.
\endaligned
\right\}
\tag 21
$$
and $\g^{(1)}_{\frac{\infty}{2}}\oplus \g^{(2)}_{\frac{\infty}{2}}$,
the highest weight vectors $w$ lying in
$Ind_k^{(l)}\setminus Ind_k^{(l+1)}$ are
$$
W^+_l = \biggl\{ \D_1^{k_1}\cdot \D_2^{k_2}\cdot \ldots \cdot v\
\vert \ \sum_{s\geq l}k_s\geq k+l,\quad \text{and}\ \sum_{s\geq
l+1}k_s\leq k+l\biggr\} .
$$
Then
$$
\chi_{\infty}(\Lambda^+_l) = \sum_{w\in
W^+_l}q^{-(k+l)l^2+k_1+4k_2+\ldots}\cdot \bigl( \chi (w)\bigr)^2.
\tag 22
$$

\head
Chapter 2.\\ The Lie Algebra \lowercase{$gl(\lambda )$}:
Irreducible Representations\\
and the Local Identity
\endhead

\def\M{\operatorname{Mat}}
\def\D{\operatorname{Det}}
\def\n{{\frak n}}
\def\p{{\frak p}}
\def\h{{\frak h}}
\def\A{{\frak A}}
\def\I{Ind_{\mu}}
\def\In{Ind_a}
\def\g{gl}
\def\e{\operatorname{deg}}
\def\s{sl}
\subhead
1. Introduction: the Lie Algebra $\g(\lambda )$ and Induced
Representations
\endsubhead
\subsubhead
1.0
\endsubsubhead
In Chapters~2 and 3 we deal with the theory of representations of the Lie
algebras $\g(\lambda ) (\lambda \in \C)$ and of the Lie algebra of
functions on a hyperboloid.  We prove the equivalence of the
irreducibility of certain representations of these Lie algebras to the
{\it local identity} (in Ch.~2) and to the {\it global identity} (in
Ch.~3) with power series (see~Subsec.~0.11). The local identity connected
with the representations of the Lie algebra $\g(\lambda )$ have the form
$$
\gather
\frac{d}{da}\Biggl( \ \prod_{i\geq 1}\frac{1}{(1-aq^i)^i}\Biggr)
\Bigg\vert_{a=1}=\\
\sum\limits\Sb{\text{over all}}\\
w=\D_1^{l_1}\cdot \ldots \cdot \D_k^{l_k}\endSb
\# \
{\Cal D}(w)\cdot q^{\sum l_i\cdot i^2}\cdot \biggl( \chi (w)\biggr)^2.
\tag 1
\endgather
$$ Here ${\Cal D}(w)$ is the Young diagram consisting of the blocks
$1\times l_1,\ldots , k\times l_k$ (see Fig.~1, Introduction), $\# {\Cal
D}(w)$ is the number of squares in it, i.e., $\sum i\cdot l_i$, and $\chi
(w)$ is a ``semiinfinite'' character corresponding to ${\Cal D}(w)$ (see
Subsec~0.11 and Ch.~1).

In Sec.~1, we introduce the principal objects of study, namely, the Lie
algebra $\g(\lambda )\ (\lambda \in \C \ )$, parabolic subalgebras in it,
and representations induced from them. The parabolic subalgebras for which
the first level of the induced representation is $k$-dimensional depend on
$k$ complex parameters and the space of possible highest weights of the
induced representation is also $k$-dimensional (see Remark~1 below).

In Sec.~2, we define the inclusion $\varphi_s:\g(\lambda )\hookrightarrow
\g_{\infty ,s}$ dependent on
$s\in \C$, where the Lie algebra $\g_{\infty,s}$ depends analytically on
$s\in \C$ and for the general $s$, is isomorphic to the Lie algebra
$\g_{\infty}$. Since $H^2(\g(\lambda ), \C \ )=0$, the inclusions
$\theta_s: \g(\lambda )\hookrightarrow
\widehat{\g}_{\infty ,s}$ are also defined. We consider the inverse image
$\theta^*_s(Ind_{\mu,s})\ (\mu \in \C \ )$ of the representation $Ind_{\mu
,s}$ of the Lie algebra $\widehat{\g}_{\infty ,s}$, which is similar to
the representation $\I$ of the Lie algebra $\widehat{\g}_{\infty}$ from
Ch.~1. We get a two-parameter family of representations of the Lie algebra
$\g(\lambda )$. We make a full investigation of the irreducibility
problem
of these representations. It turns out, on the other hand,
that in this way we get the reprersentations of $\g(\lambda )$ with a
one-dimensional first level induced from all ``maximal'' parabolic
subalgebras ${\frak p}_{\alpha}$, $\alpha \in \C$, {\it except for
two or one\/}
(according as $\lambda$).

In Sec.~3 we use the results of Ch.~1 to find the expression for the
determinant of Shapovalov's form of the representation
$\theta^*_s(Ind_{\mu ,s})$ as a function of $\mu$ and $s$. Then Theorem~2
stating that for the general $\lambda$ the representations of $\g(\lambda
)$ induced from ${\frak p}_{\alpha}$ are irreducible for one of the two
{\it exceptional\/} values of $\alpha$, and we equate, by {\it continuity}, to
zero the multiplicity of the zero of the determinant of Shapovalov's form
of the representations $\theta^*_s(Ind_{\mu ,s})$, which are close to the
exceptional representations. Thus, we get the local identity~(1).

In Sec.~4, we prove Theorem~2 on the irreducibility of representations
induced from two exceptional parabolic subalgebras for the general
$\lambda$ and $\chi (h)\ne 0$ (see also Appendix~B).
\subsubhead
1.1
\endsubsubhead
Let us consider the adjoined action of the principal $sl_2$-subalgebra of
the Lie algebra $gl_n$ on the whole algebra. Being the $sl_2$-module $$
gl_n\cong \mathop{\bigoplus}\limits^{n-1}_{i=0}\pi_i, $$ where $\pi_i$ is
an irreducible $(2i+1)$-dimensional $sl_2$-module.  Being a Lie algebra,
$gl_n$ is generated by $\pi_0, \pi_1$, and $\pi_2$, and for a sufficiently
large $n$, relations of a fixed degree depend analytically on $n$. This
allows us to determine the structure of the Lie algebra on
$\mathop{\bigoplus}\limits^{\infty}_{i=0}\pi_i$ dependent on the complex
parameter $\lambda$. The algebra $gl(\lambda )$ appearing in this way was
introduced in~[1].

On the other hand, the same Lie algebra can be constructed more
explicitly~[1].  Let us define the associative algebra $U_{\lambda}$ as a
quotient algebra $U(sl_2)/\left( \Delta - {\displaystyle\frac{\lambda
(\lambda +2)}{2}}\right)$, where $U(sl_2)$ is the universal enveloping
algebra of the Lie algebra $sl_2({\C}\ )$,\ $\lambda \in {\C}$, and
$\Delta =ef+fe+{\displaystyle\frac{h^2}{2}}\in U(sl_2)$ is the Casimir
operator.  It is easy to show that the Lie algebra (with the bracket
$[a,b]=a\cdot b-b\cdot a$) constructed using the associative algebra
$U_{\lambda}$ coincides with the Lie algebra $gl(\lambda )$.

Indeed, $\Delta$ is the central element in $U(sl_2)$ which acts on the
representation of $sl_2$ with the highest weight $\lambda$ by the
nultiplication by ${\displaystyle
\frac{\lambda (\lambda +2)}{2}}$. Thus, for
$\lambda \in {\Z}_{\geq 0}$, we have the homomorphism of Lie algebras $$
gl(\lambda )\to gl(V_{\lambda +1}), $$ where $V_{\lambda +1}$ is an
irreducible $(\lambda +1)$-dimensional $sl_2$-module with the highest
weight $\lambda$. This mapping is surjective. We denote its kernel by
$J_{\lambda}$.  Then $gl(\lambda )/J_{\lambda}\cong gl_{\lambda +1}$ and
$J_{\lambda}= \mathop{\bigoplus}\limits^{\infty}_{i=\lambda +1}\pi_i$
with respect to the adjoined action of $sl_2\subset U(sl_2)$. However, for
this definition of $gl(\lambda )$, all relations a priori depend
analytically on $\lambda$.

The Lie algebra $sl_2=\{ e,h,f\}$ is injected into $gl(\lambda )$.

We set
$$
gl(\lambda )^l=\{ v\in gl(\lambda )\ \vert \ [h,v] = 2lv\}
,\quad l\in {\Z}
$$
and
$$ {\frak
n}_{+}=\mathop{\bigoplus}\limits_{l>0}gl(\lambda )^l;\ {\frak h}=
gl(\lambda )^0;\ {\frak n}_{-}=\mathop{\bigoplus}\limits_{l<0}gl(\lambda
)^l.
$$

We have $gl(\lambda )={\frak n}_{+}\bigoplus {\frak h}\bigoplus {\frak
n}_{-}$.  Note that for $\lambda \in {\Z}_{\geq 0}$, the corresponding
decomposition of $gl_{\lambda +1} =gl(\lambda )/J_{\lambda}=
\overline{\frak n}_{+}\bigoplus
\overline{\frak h} \bigoplus \overline{\frak n}_{-}$
concides with the standard Cartan decomposition of $gl_{\lambda +1}$. In
this text, we consider representations of $gl(\lambda )$ with the highest
weight vector, i.e., a vector $v$ such that $$
\align
&{\frak n}_{+}v = 0,\\ &hv = \chi (h)\cdot v\quad
\text{for}\ h\in {\frak h}.
\endalign
$$
It is clear that the decomposition
$gl(\lambda)=\mathop{\bigoplus}\limits_{i\in {\Z}}gl(\lambda )^i$ defines the grading
on $gl(\lambda )$. Therefore, every representation with the highest weight
vector is ${\Z}_{\geq 0}$-graded.

The simplest representations with the weight vector, namely, the Verma
modules, have infinite-dimensio\-nal levels since the subalgebra ${\frak
n}_{-}$ is generated by the infinite-dimensional subset $gl(\lambda
)^{-1}$ and the Cartan algebra ${\frak h}$ acts on the levels not
semi-simply.  (The latter fact is true since $gl(\lambda )^1\quad
(gl(\lambda )^{-1})$ cannont be decomposed into the direct sum of ${\frak
h}$-invariant one-dimensional subspaces but is only entirely ${\frak
h}$-invariant.)

In order to get a discribable theory of representations, we must consider
the highest weights $\chi$ for which the corresponding Verma module has a
sufficient number of singular vectors and can be factorized  to a
representation with finite-dimensional levels.  We shall say that a
representation with finite-dimensional levels is {\it quasifinite}. These are
representations induced from parabolic subalgebras and any quasifinite
representation of $\g(\lambda)$ is a quotient of a representation induced
from some parabolic subalgebra (see Subsec.~1.3).

In addition to the $\Z$-grading $\g(\lambda )=\bigoplus\limits_{k\in
\Z}\g(\lambda )^k$, there exists a filtration
$gl(\lambda )_k=\mathop{\bigoplus}\limits^{k+1}_{i=0}\pi_i$.  The vector
space $gl(\lambda )^1\ (gl(\lambda )^{-1})$ is an analog of the space of
positive (negative) simple root vectors. It is ${\frak h}$-invariant and
cannont be decomposed into the direct sum of proper ${\frak h}$-invariant
subspaces. Note that $gl(\lambda )^{\pm 1}$ generates ${\frak n}_ {\pm}$.
\subsubhead
1.2
\endsubsubhead
Parabolic subalgebras in $\g(\lambda )$ are constructed as follows.  Let
us consider, in the ``space of negative simple root vectors'' $gl(\lambda
)^{-1}=\{ P(h)\cdot f,\ P(h)\in \C \ [h]\}$, a subspace of codimension $k$
consisting of elements of the form $P(h)\cdot f$, where $$
P=(h-\alpha_1)\ldots (h-\alpha_k)P_1,\ P_1\in {\C}\ [h].  $$ We denote by
${\widetilde {\frak n}}_-$ the Lie subalgebra in ${\frak n}_-$ generated
by this subspace.

\proclaim{Lemma~1}
${\widetilde{\frak n}}_-\bigoplus {\frak h}\bigoplus {\frak n}_+$ is a Lie
subalgebra in $gl(\lambda )$.
\endproclaim

\definition
{Definition}
The subalgebra ${\frak p}_{\alpha_1,\ldots ,\alpha_k}=
{\widetilde{\frak n}}_-\bigoplus {\frak h}\bigoplus {\frak n}_+\subset
gl(\lambda )$ is parabolic corresponding by the roots $\alpha_1,
\ldots ,\alpha_k$.
\enddefinition

Assume now that $\theta : {\frak p}_{\alpha_1,\ldots ,\alpha_k}\to {\C}$
is a one-dimensional representation. Note that this is equivalent to the
determination of the character $\theta : {\frak h}\to {\C}$ such that $$
\theta \vert_{{\frak h}\cap [\ {\widetilde{\frak n}}_-,{\frak n}_+]}=0.
\tag 2
$$
The space of these characters $\theta$ is $(k+1)$-dimensional.

Indeed, ${\frak h}\cap [{\widetilde{\frak n}}_-,{\frak n}_+]={\frak
h}\cap [ {\widetilde{\frak n}}^1_-,{\frak n}^1_+]$,
and therefore~(1)$\Leftrightarrow
\theta \vert_{[e,(h-\alpha_1)\ldots (h-\alpha_k)P_1f]}=0$
for all $P_1\in {\C}\ [h]$. It follows immediately that $\theta$
can be uniquely determined from $\theta (1),\theta (h),\ldots ,\theta (h^k)$
(see Remark~1).

Now we set
$$
L_{\alpha_1,\ldots ,\alpha_k;\theta}=U(gl(\lambda ))\mathop{\bigotimes}
\limits_{U({\frak p}_{\alpha_1,\ldots ,\alpha_k})}{\C}.
\tag 3
$$

These representations of $gl(\lambda )$ are known as
{\it generalized Verma modules}.

\proclaim{Proposition~1}
\rom{(1)}\ ${\widetilde{\frak n}}_-=\mathop{\bigoplus}\limits_{l\geq 1}
{\widetilde{\frak n}}^{\ l}_-\quad ({\widetilde{\frak n}}^{\
l}_-={\widetilde{
\frak n}}_-\cap gl(\lambda )^{-l})$, where
${\widetilde{\frak n}}_-^{\ l+1}=
\biggl\{ \biggl[ \ \prod\limits^k_{i=1}(h-\alpha_i)(h-\alpha_i+2)\ldots
(h-\alpha_i+2l)\biggr] \cdot P_1\cdot f^{l+1},\ P_1\in {\C}\
[h]\biggr\}$.

\rom{(2)} The $q$-character of the representation
$L_{\alpha_1,\ldots ,\alpha_k;\theta}$
is equal to $((1-q)(1-q^2)^2(1-q^3)^3\ldots )^{-k}$.
\endproclaim
\demo{Proof}
First, (2) follows from (1). The subspace
${\frak a}$ which is a complement of
${\widetilde{\frak n}}_-$ in
${\frak n}_-$ has $k\cdot l$ elements
of grading $-l$, and therefore, (2) follows from the
Poincar\'e--Birkhoff--Witt theorem.
To obtain (1), we note that in the associative algebra
$U(sl_2)/\biggl( \Delta -{\displaystyle \frac{\lambda (\lambda
+2)}{2}}\biggr)$, and, consequently, in $gl(\lambda )$, the relation
$f\cdot P(h)=P(h+2)\cdot f$ holds true for any $P(h)\in {\C}\ [h]$.
\enddemo
\subsubhead
1.3
\endsubsubhead
Thus, we have constructed a $(2k+1)$-parameter family of representations of
$gl(\lambda )$ with a $k$-dimensional first level. We can immediately
prove the following proposition.

\proclaim{Proposition~2}
Any quasifinite representation of $gl(\lambda )$
with a $k$-dimensional first level,
generated by one highest weight vector, is a quotient representation of
$L_{\alpha_1,\ldots ,\alpha_k;\theta}$ for certain $\alpha_1,\ldots
,\alpha_k$ and $\theta : {\frak h}\to {\C}$.
\endproclaim

\demo{Proof}
Let us consider the quasifinite representation  $N$ of
$gl(\lambda )$ with a $k$-dimensional first level, $v$ being its highest
weight vector.

We shall consider
$$
J = \{ \xi \in gl(\lambda )\ \vert \ \xi v=\mu (\xi )v,\ \mu (\xi )\in
{\C}\ \}.
$$
It is clear that $J$ is a subalgebra in $gl(\lambda ),\ J\supset {\frak
h}\bigoplus {\frak n}_+$. Let $p(h)\cdot f\in J\cap {\frak n}^1_-$.
Then $[p(h)\cdot f,\ p_1(h)]=p(h)(p_1(h+2)-p_1(h))f\in J\cap {\frak
n}^1_-$, and therefore,
$J\cap {\frak n}^1_-=\{ p(h)f,\ p\in I\}$, where $I$ is
an ideal in ${\C}\ [h]$.

It is clear that $J\cap {\frak n}^1_-\ne 0$ since the first level of $N$
is finite-dimensional and $(J\cap {\frak n}^1_-)v=0$ from considerations of
grading.

Therefore, $N$ is a quotient representation of the representation induced
from the parabolic subalgebra ${\frak p}_{\alpha_1,\ldots ,\alpha_k}$, where
$I=\{ (h-\alpha_1)\ldots (h-\alpha_k) P_1,\ P_1\in {\C}\ [h]\}$.
\enddemo
\subsubhead
1.4. Remarks
\endsubsubhead
\roster
\item  $sl(\lambda )=\mathop{\bigoplus}\limits_{i>0}\pi_i$ is a
subalgebra in $gl(\lambda )$, and we do not distinguish  between the
highest weights $\chi$ and $\chi'$ such that
$\chi \vert_{sl(\lambda )\cap {\frak h}}=\chi' \vert_{sl(\lambda )\cap
{\frak h}}$. In this sence, the representation does not depend on $\chi (1)$,
and we define a representation with a $k$-dimensional first level by
$2\cdot k$  (not by $2k+1$) parameters.

\item We can immediately show that for the general choice of  $2\cdot k$
parameters, the corresponding representation of $gl(\lambda )$ (for a fixed
$\lambda$) with the character $((1-q)(1-q^2)^2(1-q^3)^3\ldots )^{-k}$ is
irreducible. We do not prove this statement here since it follows from
the theorems given in Sec.~2.
\endroster
\subhead
2. Embeddings into $\widehat{\g}_{\infty}$ and
Critical Highest Weights
\endsubhead
\subsubhead
2.0
\endsubsubhead
In this section, we define the embeddings $\theta_s: \g(\lambda
)\hookrightarrow \widehat{\g}_{\infty ,s}$ (see Subsec.~2.1) and investigate
the problem of reducibility of the representations
$\theta^*_s(Ind_{\mu ,s})$ $(\mu ,s\in
\C \ )$.  It follows from Proposition~2 that for the given $s$, $\mu \in
\C$, there exist $\alpha$, $\chi (h)\in \C$, and a representation
$$
L_{\alpha ,\chi}\to \theta^*_s(Ind_{\mu ,s}),
$$
and both representations have the same character equal to
$\prod\limits_{i\geq 1}{\displaystyle\frac{1}{(1-q^i)^i}}$.
{\it Therefore, the irreducibility of
$L_{\alpha ,\chi}$ is equivalent to the irreducibility of
$\theta^*_s(Ind_{\mu ,s})$}, this way we obtaining all parameters
$\alpha \in \C$, except for two or one (according as $\lambda$).
\subsubhead
2.1
\endsubsubhead
The Lie algebra $gl(\lambda )$ can be constructed in a standard way with the
use of the associative algebra $U_{\lambda}=U(sl_2)/\biggl( \Delta
-{\displaystyle\frac{ \lambda (\lambda +2)}{2}}\biggr)$. This
associative algebra can be described as follows: we fix the algebra
$A={\C}\ [h]$, and being a vector space,
$$
\align
U_{\lambda}&={\C}\ [h]\bigoplus  \tag 4\\
&\bigoplus e
{\C}\ [h]\bigoplus e^2{\C}\ [h]\bigoplus \ldots \\
&\bigoplus f {\C}\ [h]\bigoplus f^2{\C}\ [h]\bigoplus \ldots.
\endalign
$$
The following relations hold true:
$$
\align
he&=e(h+2)\\
hf&=f(h-2)\\
ef&=T_1(h)=\frac{1}{2}\biggl(
h-\frac{h^2}{2}+\frac{\lambda (\lambda +2)}{2}\biggr) \\
fe&=T_2(h)=\frac{1}{2}\biggl( -h-\frac{h^2}{2}+\frac{\lambda (\lambda
+2)}{2}\biggr).
\endalign
$$

Let
$$
\align
&\sigma_1: p(h)\longmapsto p(h+2)\ \text{and}\\
&\sigma_2:
p(h)\longmapsto p(h-2)
\endalign
$$
be two automorphisms of the algebra $A={\C}\ [h]$, and then
$$
\left\{
\aligned
&p\cdot e=e\cdot \sigma_1(p)\\
&p\cdot f=f\cdot \sigma_2(p)\\
&ef=T_1;\ fe=T_2\\
&\sigma_1T_1=T_2
\endaligned
\right.
\tag 5
$$
where $p,T_1,T_2\in {\C}\ [h]$.

Let us now consider the associative algebra $Mat_{\infty}$ of generalized
Jacobian matrices. We set\newline
$
A=\prod\limits^{\infty}_{i=-\infty}{\C}\ (i),
$
where ${\C}\ (i)$ is a copy
of the algebra ${\C}$. Being a vector space,
$$
\align
Mat_{\infty}&=A\bigoplus \tag 6\\
&\bigoplus eA\bigoplus
e^2A\bigoplus e^3A\bigoplus \ldots \\
&\bigoplus fA\bigoplus
f^2A\bigoplus f^3A\bigoplus \ldots
\endalign
$$
$$
\left\{
\aligned
&ef=fe=1\in A\\
&H\cdot e=e\cdot D(H),\ \text{where}\\
&D: A\to A\ -\
\text{a shift to the right by unity}.
\endaligned
\right.
\tag 7
$$

(Here $A$ are diagonal matrices, $e$ is a diagonal of unities which is
above the principal, and $f$ is a diagonal of unities which is below the
principal diagonal.)

The main difference between $Mat_{\infty}$ and $U_{\lambda}$ is that
$T_1=T_2=1$ in (7), and then, as in (5), $T_i\ne 1$.

It follows that we can construct the mapping of the associative algebras
$\varphi_s: U_{\lambda}\to Mat_{\infty ,s}$, where, being a vector space,
$$
\align
Mat_{\infty ,s}&=A\bigoplus \tag 8\\
&\bigoplus eA\bigoplus
e^2A\bigoplus e^3A\bigoplus \ldots \\
&\bigoplus fA\bigoplus
f^2A\bigoplus f^3A\bigoplus \ldots
\endalign
$$
where $A=\prod\limits^{\infty}_{i=-\infty}{\C}\ (i)$ and
$$
\gather
\left\{
{\aligned
&H\cdot e=e\cdot D(H)\\
&H\cdot f=f\cdot D^{-1}(H)\\
&ef=(T_1(s-2i))_{i\in {\Z}}\in A\\
&fe=(T_2(s-2i))_{i\in {\Z}}\in A,
\endaligned}
\right.
\tag 9\\
s\in \C
\endgather
$$
$Mat_{\infty ,s}$ is an associative algebra by virtue of the relation
$\sigma_1 T_1=T_2$. In order to construct this mapping,
we must construct the mapping $\varphi_s: {\C}\ [h]\to
\prod\limits^{\infty}_{i=-\infty}{\C}\ (i)$
such that the diagram
$$
\CD
{\C}\ [h]  @>{\varphi_s} >>
\prod\limits^{\infty}_{i=-\infty}{\C}\ (i)\\
@V\sigma_1VV @VVDV\\
{\C}\ [h]  @>{\varphi_s}>> \prod\limits^{\infty}_{i=-\infty}{\C}\ (i)
\endCD
$$
is commutative.

It follows that if $\varphi_s(h)_0=s\in {\C}\ (0)$, then
$\varphi_s(h)_i= s-2i\in {\C}\ (i)$. Thus
$$
\varphi_s(p(h))_i=p(s-2i)\in {\C}\ (i),
\tag 10
$$
where $p\in {\C}\ [h]$.

\proclaim
{Lemma~2}
$\varphi_s: U_{\lambda}\to Mat_{\infty ,s}$ is an embedding.
\endproclaim

We have defined the embedding of the corresponding Lie algebras
$\varphi_s: gl(\lambda)\hookrightarrow gl_{\infty ,s}$.
We will define the associative algebra $Mat_{\infty ,(i)}\ (i\in {\Z})$
by relations~(8) and (9), where, instead of the relations for
$ef$  and $fe$ in (9), we set for $j\ne i$ and $(ef)_i=0$;
$ef\in A$,\ $(ef)_j=1$ for
$fe\in A$,\ $(fe)_j=1$ for $j\ne i+1$ and
$(fe)_{i+1}=0$. Clearly, the algebra $Mat_{\infty ,(i)}$ is associative, and
we denote by $gl_{\infty ,(i)}$ the corresponding Lie algebra.

In what follows, we distinguish between two cases, namely,
$$
\left\{
\alignedat 2
&\text{(i)}\quad&T_1(s-2i)=0\ &\text{ for a certain }\ i\in {\Z},\\
&\text{(ii)}\quad&
T_1(s-2i)\ne 0\ &\text{ for all }\ i\in {\Z}.
\endalignedat
\right.
\tag 11
$$

\proclaim{Lemma~3} For the general $\lambda$,

$Mat_{\infty ,s}\cong Mat_{\infty, (i)}$ in case \rom{(i)},

$Mat_{\infty ,s}\cong Mat_{\infty}$ in case \rom{(ii)}.
\endproclaim

\demo{Proof}
The map $\Psi_s: Mat_{\infty ,s}\to Mat_{\infty}$ is defined:
$$
\aligned
H&\longmapsto H \\
f&\longmapsto f\\
e&\longmapsto (T_1(s-2i),\ i\in {\Z})\cdot e,
\endaligned
\tag 12
$$
which is an isomorphism in case~(ii). This proves~(ii),
and case~(i) can be proved by analogy.
\enddemo
The assumption concerning the genaral $\lambda$ is that $T_1$ does not
have two roots $t_1$ and $t_2$ such that $t_1-t_2\in 2{\Z}$.
When these two roots exist, we have
$Mat_{\infty ,s}\cong Mat_{\infty ,(i,j)}$ for certain $s\in {\C}$, where
the last associative algebra can be determined in an obvious way.

\definition{Definition}
The associative algebra $U^{\Cal O}_{\lambda}$
can be constructed from relations~(4) and (5) where $A={\C}\ [h]$ is replaced
by $A={\Cal O}$, the algebra of holomorphic functions
${\C}\to {\C}$. The Lie algebra
$gl^{\Cal O}(\lambda )$ is the Lie algebra corresponding to the
associative algebra $U^{\Cal O}_{\lambda}$.
\enddefinition

\proclaim{Lemma~4}
The embedding $\varphi_s: gl(\lambda )\hookrightarrow
gl_{\infty ,s}$ can be continued to the map $\varphi_s^{\Cal O}: gl^{\Cal
O}(\lambda )\to gl_{\infty ,s}$ which is surjective.
\endproclaim

\demo{Proof}
The first statement is obvious and the second statement follows
from the fact that
for the discrete sequence $(s-2i)_{i\in {\Z}}$ there exists a
holomorphic function from
$\C$ into $\C$ which assumes the given values at the points of this
sequence.
\enddemo
\subsubhead
2.2
\endsubsubhead
There is a 2-cocycle $\alpha$ of the Lie algebra $gl_{\infty ,s}$
defined as
$$
\aligned
\alpha (E_{ij},E_{ji})&={\cases
P_{j-i}(s-2i)\ &\text{ for }\ i\leq 0, j\geq 1\\
0\ &\text{ otherwise }
\endcases}\\
\alpha (E_{ij}, E_{kl})&=0\ \text{for}\ i\ne l\ \text{and}\ j\ne k.
\endaligned
\tag 13
$$

Here $P_l(h)=e^lf^l,\ E_{ij}=1_ie^{j-i}$ for $i<j$, where $1_i\in
\prod\limits^{\infty}_
{i=-\infty}{\C}\ (i)$ is a sequence with 1 is the $i$th position and 0
in the other positions. Note that in case~(ii) (see~(11)) the cocycle $\alpha$
is the inverse image for $\Psi_s$ (see~(12)) of a standard cocycle on
$gl_{\infty}$ and in case~(i) the cocycle $\alpha$ is coholomogous to zero.

The $\widehat{\g}_{\infty ,s}$-module  $Ind_{\mu ,s}$ can be termined
in a standard way as a quotient module of the Verma module
$M_{\mu}$ with the zero highest weight and the central charge
$\mu \in \C$, by analogy with the module of $\I$ in the
$\widehat{\g}_{\infty}$-case (see Ch.~1).

\proclaim{Lemma~5}
\roster
\item If $T_1(s-2i)=0$ for a certain $i\in \Z$ \rom(case \rom{(i))}, then
$Ind_{\mu ,s}$ is reducible for all $\mu \in \C$,
\item If $T_1(s-2i)\ne 0$ for all $i\in \Z$ \rom(case \rom{(ii))} and $\mu
\in \Z$, then $Ind_{\mu ,s}$ is reducible,
\item If $T_1(s-2i)\ne 0$ for all $i\in \Z$ \rom(case \rom{(ii))} and $\mu
\notin \Z$, then $Ind_{\mu ,s}$ is irreducible.
\endroster
\endproclaim

\demo{Proof}
According to Lemma~3, statement~(1) reduces to the corresponding statement
concerning $\g_{\infty ,(i)}$, and (2) and (3) reduce to the corresponding
statements concerning $\widehat{\g}_{\infty}$. Cases (2) and (3) follow
directly from Theorem~1, Ch.~1; in general, all the statements follow from the
results of Sec.~3 of this Chapter, where we calculate the determinant
of Shapovalov's from of the representation
$\theta^*_s(Ind_{\mu ,s})$. However, we give here the direct proof of~(1)
for convenience.

We set $e_l=1_l\cdot e;\ f_l=1_{l+1}\cdot f$. These are the positive and the
negative root vectors in $gl_{\infty ,(i)}$. For $gl_{\infty,(i)}$ we have
$$
[e_i, f_i] = 0.
\tag 14
$$
We set $x=[\ldots [[f_0,f_1],\ f_2],\ldots ,f_i]$ for $i\geq 0$ and
$x=[\ldots [[f_i,f_{i+1}],\ f_{i+2}],\ldots ,f_0]$ for $i\leq 0$
\vskip+5cm

\centerline{Fig.~1}

In Fig.~1\ $x$ is a hatched element. Let $v$ be the highest weight
vector in $Ind_{\mu ,s}$. Then $xv$ is a singular vector in $Ind_{\mu
,s}$.  Indeed, it is sufficient to verify that
$$
e_0xv=0\ \text{in}\ Ind_{\mu ,s}
\tag 15
$$
and
$$
e_ixv=0\ \text{in}\ Ind_{\mu ,s},
\tag 16
$$
since if $j\ne 0,\,i$ then $[e_j,x]=0$ in $gl_{\infty ,(i)}$.
In order to prove~(15), we note that $[e_0,x]$
lies outside of the angle hatched in Fig.~2,
and therefore, $[e_0,x]v=0$ in $Ind_{\mu ,s}$; (16) follows from the
permutability of $e_i$ with all $f_j$,\ $j\in {\Z}$ in $U(gl_{\infty,(i)})$.
\enddemo
\subsubhead
2.3
\endsubsubhead
We can show that $H^2(gl(\lambda ),{\C}\ )=0$, i.e.,
$gl(\lambda )$ does not have any nontrivial central extensions. Therefore,
the embedding
${\widehat{\varphi}}_ s:{\widehat{gl}}(\lambda )\hookrightarrow
{\widehat{gl}}_{\infty ,s}$ with an induced cocycle on $gl(\lambda )$
defines the embedding $\theta_s: gl(\lambda
)\hookrightarrow{\widehat{gl}}_{\infty ,s}$. From~(13) it is easy to find
$$
\left\{
\aligned
&\theta_s(h)=\varphi_s(h)+T_1(s)\cdot c\\
&\theta_s(3h^2-2c_{\lambda})=
\varphi_s(3h^2-2c_{\lambda})+(2s-2)T_1(s)\cdot c,
\endaligned\right.
\tag 17
$$
where $c_{\lambda}={\displaystyle\frac{\lambda (\lambda +2)}{2}}$,
$T_1(s)= {\displaystyle\frac{1}{2}}\biggl(
h-{\displaystyle\frac{h^2}{2}}+c_{\lambda}\biggr)$, and $c$ is the central in
${\widehat{gl}}_{\infty ,s}$. On the other hand, $\theta
(1)=1+\varepsilon \cdot c$, and in a nondegenerate case, the embeddings
$\theta_s$ are parametrized by the choice of $\varepsilon \in {\C}$.
There is the commutative diagram
$$
\CD
{\widehat{gl}}^{\Cal O}(\lambda ) @>
\varphi^{\Cal O}_s>>  {\widehat{gl}}_{\infty ,s}\\
@AAtA @.\\
gl(\lambda ) @.\\
\endCD
$$
where $t$ is a natural embedding which preserves the grading and
${\widehat{\varphi}}^ {\Cal O}_s$ is surjective. We shall denote the
grading by a superscript.

\proclaim{Lemma~6} Let $V$ be a quasifinite $gl(\lambda
)$-module.  Then the action of $gl(\lambda )^k$ on $V$ for $k\ne 0$
can be naturally continued to the action of
${\widehat{gl}}^{\Cal O}(\lambda )^k$
on $V$.
\endproclaim

\proclaim{Corollary} Any $gl(\lambda )$-invariant subspace $W$
of the module $Ind_{\mu ,s}$ is invariant with respect to
$\widehat{gl}_{\infty ,s}$.
\endproclaim

\demo{Proof of the Corollary}
Since ${\widehat{\varphi}}^{\Cal O}_s: {\widehat{gl}}^{\Cal O}(\lambda
)^k
\to {\widehat{gl}}^k_{\infty ,s}$ is surjective, it follows, in
accordance with Lemma~6 that $W$ is also
${\widehat{gl}}^k_{\infty ,s}$-invariant for $k\ne 0$, whence follows the
${\widehat{gl}}_{\infty ,s}$-invariance.
\enddemo

\demo{Proof of the Lemma~\rm[2]} Let $V$ be a graded quasifinite
representation:
$V=\mathop{\bigoplus}
\limits_{p\geq 0}V_p$, $dim\ V_p<\infty$. We shall consider
$Hom\ (V,V):=\mathop{
\bigoplus}\limits_{p,q\geq 0}Hom\ (V_p,V_q)$ with the topology of the
direct sum of finite-dimensuonal spaces. We can assume that
$V_p$ are normed spaces, and then  the norm $\lnorm ,\rnorm_{p,q}$
is induced on $Hom\ (V_p,V_q)$. We shall show that for $k\ne 0$,
the mapping  $gl(\lambda )^k\to Hom\ (V,V)$
is continuous; for which purpose we shall show that
$\lnorm e^kh^n\rnorm$ is bounded in
$Hom\ (V_p,V_{p+k})$ for fixed $k$, $p$, and an arbitrary $n$. We note
that $e^kh^n={\displaystyle
\frac{1}{(4k)^n}}(adh^2-4k^2)^ne^k$, where $e^k\in Hom\ (V_p,V_{p+k})$, and
$(adh^2-k^2)\in End\ (Hom\ (V_p,V_{p+k}))$.

We have  $\lnorm e^kh^n\rnorm_{p,p+k}\leq A\cdot \alpha^n$, where $A=\lnorm
e^k\rnorm$ and
$\alpha ={\displaystyle\frac{\lnorm adh^2-k^2\rnorm}{2k}}$,
and therefore $\lnorm e^kf(h)\rnorm_{p,p+k}=\lnorm \sum\limits_{n\geq
0}f_ne^kh^n\rnorm_ {p,p+k}\leq \sum\limits_{n\geq 0}\lmod f_n\rmod \cdot
\lnorm e^kh^n\rnorm_
{p,p+k}\leq A\cdot\sum\limits_{n\geq 0}\lmod f_n\rmod
\cdot \alpha^n=A\cdot \varphi (f)(\alpha )$.

The statement of the Lemma follows from the fact that
$\varphi : \sum f_iz^i\longmapsto
\sum \lmod f_i\rmod z^i$ is continuous and that $\Cal O$
is a completion of ${\C}\ [z]$ in the topology of the uniform convergence
on compact sets.
\enddemo
\subsubhead
2.4
\endsubsubhead
We have proved that the  $\g(\lambda )$-module $\theta^*_s(Ind_{\mu
,s})$ is reducible for $\mu \in \Z$ or for
$T_1(s-2i)=0$ for a certain $i\in
\Z$ and irreducible otherwise. Its highest weight $\chi$
is
$$
\left\{
\aligned
&\chi (1)=\varepsilon \cdot \mu \\
&\chi (h) = T_1(s)\cdot \mu \\
&\chi (3h^2-2c_{\lambda})=(2s-2)T_1(s)\cdot \mu.
\endaligned
\right.
\tag 18
$$

Let us compare this with the expression for the representation
induced from ${\frak p}_{\alpha}$
with the given $\chi (1).\ \chi (h)$:
$$
\left\{
\aligned
&\chi (1)= \chi (1)\\
&\chi (h) = \chi (h) \\
&\chi (3h^2-2c_{\lambda})=2(1+\alpha )\cdot \chi (h).
\endaligned
\right.
\tag 19
$$

We compare~(12) with (13) and express $(\alpha ,\chi (h))$ parameters
in terms of the $(s,\mu )$ parameters.
If $T_1(s)\ne 0$; then
$$
s=\alpha +2\ \text{and}\ \mu =\chi (h)\big/ T_1(\alpha +2).
\tag 20
$$
It is now clear that{\it if\ $T_1(\alpha +2)=0$,\ then,\ for\ all\ $\chi
(h)$,\ the corresponding representation of\ $\g(\lambda )$ does not have
the $(s,\mu)$-parametrization}.

We have proved the following theorem
\proclaim{Theorem~1} For $T_1(\alpha +2)\ne 0$ the representation of the Lie
algebra  $\g(\lambda )$, induced from the parabolic subalgebra
${\frak p}_{\alpha}$ with the highest weight $\chi$, is reducible if

\rom {(1)} $T_1(\alpha +2i)=0$ for the integer $i\ne 1$\newline
or

\rom {(2)} for ${\displaystyle\frac{\chi (h)}{T_1(\alpha +2)}}\in\Z$.\newline
In other cases, this representation is irreducible.
\endproclaim

\proclaim{Theorem~2}
For all $\chi (h)\ne 0$, the representation of the Lie algebra
$\g(\lambda )$ induced from ${\frak p}_{\alpha}$ is reducible for
$T_1(\alpha +2)=0$
if $T_1(\alpha +2i)=0$ for a certain integer $i\ne 1$ and irreducible
otherwise.
\endproclaim

As we shall see in Sec.~3, the local identity~(1) is equivalent to this
theorem.
The proof of the Theorem~2 will be given in Sec.~4.
\remark
{2.5. Remark}
Using the technique developed, we can easily show that the
$\g(\lambda )$-module $\bigotimes\limits^k_{i=1}L_{\alpha_i,\chi_i}$
{\it does not\/} coincide with the generalized Verma module with the same
highest weight if $\alpha_i-\alpha_j\in 2\cdot \Z$ for certain $i\ne j$ and
coincides with it otherwise. This module is reducible if
$L_{\alpha_i,\chi_i}$ is reducible for a certain $i$ or if
$\alpha_i-\alpha_j\in 2\cdot \Z$ for certain $i\ne j$ and irreducible
otherwise.
\endremark
\subhead
3. The Local Identity
\endsubhead
\subsubhead
3.0
\endsubsubhead
We fix  a general (see Theorem~2 in Subsec.~2.4) $\lambda \in
\C$, and let $\alpha$ be such that $T_1(\alpha +2)=0$. It is clear
that the determinant of Shapovalov's form depends analytically on
$\alpha$ and $\chi (h)$
(in the sense refined below this determinant is not uniquely defined).
Therefore, Theorem~2 from Subsec.~2.4 is equivalent to the statement that for
$T_1(s)\sim t$ and $\mu \sim {\displaystyle\frac{1}{t}}$ as $t\to 0$ (see
(20)), the corresponding limit of the determinant of Shapovalov's form
does not depend on  $t$, i.e., tends to a certain finite nonzero number.
\subsubhead
3.1. Chevalley involution
\endsubsubhead
Let $\omega$ be the Chevalley involution of the Lie algebra
$\s_2: \omega (e)=-f,\ \omega (f)=-e,\ \omega (h)=-h$. There exist two
continuations of $\omega$ to ${\Cal U}(\s_2)$ as a Lie algebra, namely,
$$
\alignat 2
&\text{(i)}&\quad & \omega_1(1)=1,\ \omega_1\big\vert_{\s_2}=\omega ,\\
&&\quad&\omega_1(a\cdot b)=\omega_1(a)\cdot \omega_1(b);\\
&\text{(ii)}&\quad& \omega_2(1)=-1,\ \omega_2\big\vert_{\s_2}=\omega ,\\
&&\quad&\omega_2(a\cdot b)=-\omega_2(b)\cdot \omega_2(a).
\endalignat
$$
For our purposes (i.e., for the existence of a map from the Verma module
into a countergradient module) the condition $\omega
\big\vert_{\h}=-Id$ is necessary, which is satisfied for
$\omega_2$. It is also clear that
$\omega_2$ defines the involution of the Lie algebra $\g(\lambda )$ which
is known as the Chevalley involution, and is denoted by $\omega$.
\subsubhead
3.2
\endsubsubhead
Knowing the expression for the determinsant of Shapovalov's form
of the $\widehat{\g}_{\infty}$-module of $\I$ as a function of
$\mu \in \C$, we want to find the expression for the determinant of
Shapovalov's form of the $\g(\lambda )$-module
$\theta^*_s(Ind_{\mu ,s})$ as a function of $\mu ,s\in \C$.

We denote $\Gamma_s=\Psi_s\varphi_s$ (see (10), (12)),
$\Gamma_s: U_{\lambda}\hookrightarrow Mat_{\infty}\ (s\in \C \ )$.
Recall that $\Gamma_s$ is defined as follows:
$$
\left.
\aligned
&e\mapsto \biggl( T_1(s-2i),\ i\in \Z \biggr) \cdot e\\
&p(h)\mapsto \biggl(
p(s-2i),\ i\in \Z \biggr) \\
&f\mapsto f.
\endaligned
\right\}
\tag 21
$$

Assuming that $s\in \C$, we define the holomorphic functions
$\delta_{i,s}$ on $\C$  by the condition
$$
\left\{
\aligned
&\delta_{i,s}(s-2i)=1,\\
&\delta_{i,s}(s-2j)=0\ \text{for}\ j\ne i.
\endaligned
\right.
$$
(Clearly, these functions are not uniquely defined.)

For our purposes $\{ f^k\delta_{l,s}\} \subset {\n}_-^{(-k)}\subset
\g^{\Cal O}(\lambda )$ is more convenient than the standard base
$\{ f^kh^l\}$ since under the surjection
$\Gamma_s^{\Cal O}: \g^{\Cal O}(\lambda )\to \g_{\infty}$
the elements of the first set pass into the standard base $\{ E_{ij},\
i>j\} \subset \n_-\subset \g_{\infty}$ (and this is very convenient for
calculating of Shapovalov's form). {\it In this case, the multiplicities of
the roots of the determinant of Shapovalov's form do not change
when $\g(\lambda )$ is replaced by $\g^{\Cal O}(\lambda )$}.
\subsubhead
3.3
\endsubsubhead
It is clear that the base on the $k$th level of the induced representation
of $\g(\lambda )$ with a one-dimensional first level is formed by the vectors
$$
w=\bigl( f^{i_1}\delta_{j_1,s}\bigr) \cdot \ldots \cdot \bigl(
f^{i_p}\delta_{j_p,s}\bigr) \cdot v,
\tag 22
$$
where $v$ is the highest weight vector, $i_1+\ldots +i_p=k$, and
$$
j_q\in [-i_q+1,0]\ \text{for all}\ q.
\tag 23
$$
We want to calculate Shapovalov's form in this base. It is clear that
$$
\Phi_{\g^{\Cal O}(\lambda )}(w_1,w_2)=\Phi_{\widehat{\g}_{\infty}}\bigl(
\Gamma^{\Cal O}_sw_1, \Gamma^{\Cal O}_sw_2\bigr) ,
\tag 24
$$
where the $\g(\lambda )$-module on the left-hand side is the inverse image
of the $\widehat{\g}_{\infty}$-module for $\Gamma_s$ on the right-hand side.
In explicit form,
$$
\gather
\Phi_{\g^{\Cal O}(\lambda )}(w_1,w_2)=\\
\pm \bigl( \delta_{j'_q,s}e^{i'_q}\bigr) \cdot \ldots \cdot
\bigl( \delta_{j'_1,s}e^{i'_1}\bigr) \cdot \bigl(
f^{i_1}\delta_{j_1,s}\bigr) \cdot \ldots \cdot \bigl(
f^{i_p}\delta_{j_p,s}\bigr) v\\
\text{(for\ the representation of}\ \g(\lambda ))=\\
\pm \Gamma^{\Cal O}_s\bigl( \delta_{j'_q,s}e^{i'_q}\bigr) \cdot \ldots
\cdot \Gamma^{\Cal O}_s\bigl( \delta_{j'_1,s}e^{i'_1}\bigr) \cdot
\Gamma^{\Cal O}_s\bigl( f^{i_1}\delta_{j_1,s}\bigr) \cdot \ldots \cdot
\Gamma^{\Cal O}_s\bigl( f^{i_p}\delta_{j_p,s}\bigr) v\tag 25\\
\text{(for\ the representation of}\ \widehat{\g}(\lambda )).
\endgather
$$
We have
$$
\gather
\Gamma^{\Cal O}_s\bigl(\delta_{j,s}e^i\bigr) =\Gamma^{\Cal O}_s\bigl(
\delta_{j,s}\bigr) \cdot \Gamma^{\Cal O}_s(e)\cdot \ldots \cdot
\Gamma^{\Cal O}_s(e)=\\
\BI_{j}\cdot (T)\cdot e\cdot (T)\cdot e\cdot \ldots \cdot (T)\cdot e,
\tag 26
\endgather
$$
where $\BI_{j}\in \prod\limits_{i\in \Z}\C \ (i)$ is a sequence with 1 in
the  $j$th place and with 0 in the other places and $(T)=\biggl(
T_1(s-2i),\ i\in \Z \biggr) \in \prod\limits_{i\in \Z}\C
\ (i)$. On the other hand,
$$
\Gamma^{\Cal O}_s \bigl( f^i\delta_{j,s}\bigr) =f^i\cdot \ \BI_{\ j}.
\tag 27
$$
It is obvious from~(26) that
$$
\Gamma^{\Cal O}_s\bigl(
\delta_{j,s}e^i\bigr) =\biggl( \ \prod^{i-1}_{l=0}T_1(s-2j-2l)\biggr)
\cdot \ \BI_{\ j}\cdot e^i,
\tag 28
$$
and the determinant of Shapovalov's form of the $\g(\lambda )$-module
$\theta^*_s(Ind_{\mu ,s})$ is
$$
\gather
\det \Phi_{\g(\lambda )}(s,\mu )=\\
\Biggl( \ \prod\Sb\text{over all}\\
w\ \text{from the}\ k\text{th}\\
\text{level~(22)}\endSb\prod^p_{m=1}\ \prod^{i_p-1}_{l=0}
T_1(s-2j_m-2l)\Biggr) \cdot \det \Phi_{\widehat{\g}_{\infty}}(\mu ).
\tag 29
\endgather
$$
Recall (see~(23)) that $j_q\in [-i_q+1,0]$.
\subsubhead
3.4
\endsubsubhead
Assume now that $t\to 0$, and let $T_1(s)\sim t$ and $\mu \sim
{\displaystyle\frac{1}{t}}$ (see Subsec.~3.0). Then $\det
\Phi_{\widehat{\g}_{\infty}}(\mu )\sim
{\displaystyle\frac{1}{t^{\deg_k}}}$, where $\deg_k$ is the degree of $\det
\Phi_{\widehat{\g}_{\infty}}(\mu )$ for $\mu$ on the $k$th level.
On the other hand, according to condition~(23), the multiplicity of zero
of the product in relation~(29) for a general $\lambda$ is equal to
$p$ for every $w=\bigl(
f^{i_1}\delta_{j_1,s}\bigr) \cdot \ldots \cdot \bigl(
f^{i_p}\delta_{j_p,s}\bigr)$ since $T_1(s-2i)$ has 0 only for $i=0$.

Thus, according to Lemma~12 from Subsec.~2.2 of Ch.~1, Theorem~2
from Subsec.~2.4 is equivalent to the local identity~(1).
\subhead
4. Proof of the Local Identity: Embeddings of
$\g(\lambda )$ into $\widehat{\g}_{\infty}(\C \
[t]/t^2)$
\endsubhead
\subsubhead
4.0
\endsubsubhead
In Sec.~3, we have derived the local identity~(1) from Theorem~2. This section
is devoted to the proof of the theorem itself. First, Theorem~2 follows
from the following weaker result.
\proclaim{Theorem~3}
For a general $\lambda$, $T_1(\alpha +2)=0$, and
$\chi (h)\ne 0$, the representation of the Lie algebra
$\g(\lambda )$ induced from ${\frak p}_{\alpha}$ is irreducible.
\endproclaim

Indeed, the local identity follows from Theorem~3 by analogy with Sec.~3.
However, the conditions of Theorem~2 are precisely the conditions which
are necessary for the arguments from Sec.~3, in order to have no additional
zeros.

In this section, we prove Theorem~3. As we have seen in Sec.~2, for
$T_1(\alpha +2)=0$, the representation of the Lie algebra $\g(\lambda )$
induced from ${\frak p}_{\alpha}$ does not have the
$(s,\mu )$-parametrization.

However, we can obtain its highest weight by considering the embeddings
$\theta_s: \g(\lambda )\hookrightarrow {\widehat \g}_{\infty ,s}(\C
\ [t]/t^2)$ (see below) and the inverse images of the
${\widehat \g}_{\infty ,s}(\C
\ [t]/t^2)$-modules  $Ind_{\mu_1,\mu_t,s}$ for $\theta_s;\
\theta^*_s(Ind_{\mu_1,\mu_t,s})$. (Here $\mu_1,\mu_t$ are two central
charges.) For certain $s$, $\mu_1$, $\mu_t$, this representation
has the same highest weight as that induced from the exceptional parabolic
subalgebra ${\frak p}_{\alpha}$, and we prove that the irreducible
quotient module for these $\mu_1$, $\mu_t$, $s$ of the module
$\theta^*_s(Ind_{\mu_1,\mu_t,s})$ has, for
$\mu_t\ne 0$, a character not smaller than $\prod\limits_{i=1\ldots
\infty}(1-q^i)^{-i}$. On the other hand, this character is not greater
than $\prod\limits_
{i=1\ldots \infty}(1-q^i)^{-i}$ since there exists a representation
induced from ${\frak p}_{\alpha}$ with the corresponding highest weight.
Hence Theorem~3 follows.
\subsubhead
4.1
\endsubsubhead
We shall consider the ``finite'' algebra $\g_{2n} (\C \
[t]/t^2)$ instead of $\widehat{\g}_{\infty}(\C \ [t]/t^2)$ and then pass
to the limit. We denote by
$\overline{\g}_{2n}:=\g_{2n}\!\cdot\! t$,
$\overline{\g}_{2n}$ the Abelian Lie algebra, being a $\g_{2n}$-module,
$\overline{\g}_{2n}$ is isomorphic to the adjoint action on $\g_{2n}$.
There are two subalgebras, $\g^{(1)}_n$ and $\g^{(2)}_n$, in
$\g_{2n}$, and they commute with each other, and there are also two
Abelian subalgebras  $\A_+$ and $\A_-$ (see Fig.~1 in Ch.~1).
We assume a similar notation for $\overline{\g}_{2n}$,
$\g_{2n}=\n_-\oplus \h\oplus \n_+$,
$\g^{(1)}_n=\n_-^{(1)}\oplus \h^{(1)}\oplus \n_+^{(1)}$ and so on.
\subsubhead
4.1.1
\endsubsubhead
The associative algebra $\M_{\infty}(\C \ [t]/t^2)$ is described as follows.

We denote by $A=\prod\limits_{i\in \Z}(\C \ [t]/t^2)_i$ the direct product of
algebras. Being a vector space,
$$
\align
\M_{\infty}(\C \ [t]/t^2)=& A\oplus \\
& \oplus eA\oplus e^2A\oplus e^3A\oplus \ldots \\
& \oplus fA\oplus f^2A\oplus f^3A\oplus \ldots,
\endalign
$$
and the relations
$$
\left\{
\aligned
&H\cdot e=e\cdot D(H)\\
&H\cdot f=f\cdot D^{-1}(H)\\
&e\cdot f=(\ldots ,1,1,1,\ldots )\in A\\
&f\cdot e=(\ldots ,1,1,1,\ldots )\in A,
\endaligned
\right.
\tag 30
$$
are satisfied, where $D: A\to A$ is a shift to the right by~1.

We denote by $a_i\in A$ a sequence with $t$ in the  $i$th place and~1
in the other places. Then, replacing the last two conditions
in~(30) by
$$
\left\{
\aligned
&.\ .\ .\ .\ .\ .\ .\ .\ .\ .\\
&e\cdot f=a_i\in A\\
&f\cdot e=a_{i+1}\in A,
\endaligned
\right.
\tag 31
$$
we get the associative algebra $\M_{\infty ,i}(\C \ [t]/t^2)$.
We denote by $\g_{\infty ,i}(\C \ [t]/t^2)$ the corresponding Lie algebra.
It is clear that there exists a ``finite'' analog of $\g_{2n,i}(\C
\ [t]/t^2)$.
\subsubhead
4.1.2
\endsubsubhead
Let $\alpha_0^{\vee},\ \overline{\alpha_0^{\vee}}$ be the
``central'' coroots in $\g_{2n}(\C \ [t]/t^2)$. We shall considers
the $\g_{2n}(\C\ [t]/t^2)$-module of $Ind_{\mu_1,\mu_t,s}$ induced from
$\p\oplus \overline{\p}$ with the highest weight $\chi$ such that $\chi
(\alpha_0^ {\vee})=\mu_1$, $\chi (\overline{\alpha_0^{\vee}})=\mu_t$.
Here $\p=\g_n^{(1)}\oplus \A_+\oplus \g_n^{(2)}$ and
$\overline{\p}=\overline{\g_n^{(1)}}\oplus \overline{\A_+}\oplus
\overline{\g_n^{(2)}}$. We denote by $Ind_{\mu_1,\mu_t,s}$ a similar
module over $\g_{2n,i}(\C \ [t]/t^2)$
($\widehat{\g}_{\infty ,i}(\C \ [t]/t^2)$, respectively).

\proclaim{Lemma~7}
The $q$-character of the ireducible quotient module of the
$\widehat{\g}_{\infty ,0}(\C \ [t]/t^2)$-module  $Ind_{\mu_1,\mu_t,s}$
is greater or equal to
$\prod\limits_{i\geq 1}(1-q^i)^{-i}$ for $\mu_t\ne 0$.
\endproclaim

\remark{Remark}
We shall see in Subsec.~4.2.3 that for $\mu_t\ne 0$ this character
is equal to $\prod\limits_{i\geq 1}(1-q^i)^{-i}$.
\endremark
\demo
{4.1.3. Proof}
Let us consider the action
$\g^{(1)}_{\frac{\infty}{2}}\oplus \g^{(2)}_{\frac{\infty}{2}}\subset
\widehat{\g}_{\infty ,0}(\C \ [t]/t^2)$ on the module
$Ind_{\mu_1,\mu_t,s}$. Being a vector space,
$Ind_{\mu_1,\mu_t,s}\cong S^*(\A_-)\otimes S^*(\overline{\A_-})$.
According to what was stated in Ch.~1 (Lemma~6 in Sec.~2), the
$\g^{(1)}_{\frac{\infty}{2}}\oplus
\g^{(2)}_{\frac{\infty}{2}}$-module $S^*(\A_-)$ is isomorphic
to the direct sum of all irreducible modules $L_w\otimes L_w$
generated by the highest weight vectors
$w=\D^{k_1}_1\cdot \ldots \cdot \D^{k_l}_l\cdot v$
(here $v$ is the highest weight vector in $Ind_{\mu_1,\mu_t,s}$ and $\D_k\in
S^*(\A_-)$ is the determinant of the matrix
$A_k=(y_{ij},\ i,j=1\ldots k)$).
Therefore, $Ind_{\mu_1,\mu_t,s}\cong S^*(\A_-)\otimes
S^*(\overline{\A_-})$ is isomorphic to the direct sum of the tensor
products of irreducible $\g^{(1)}_{\frac{\infty}{2}}\oplus
\g^{(2)}_{\frac{\infty}{2}}$-modules $(L_w\otimes L_w)\otimes
(L_{\overline w}\otimes L_{\overline w})$, each of which, in turn,
can be decomposed into the direct sum of irreducible
$\g^{(1)}_{\frac{\infty}{2}}\oplus \g^{(2)}_{\frac{\infty}{2}}$-modules.
We denote by $w\cdot {\overline w}\cdot v$ the highest weight vector of
$(L_w\otimes
L_w)\otimes (L_{\overline w}\otimes L_{\overline w})$.
The $\g^{(1)}_{\frac{\infty}{2}}\oplus
\g^{(2)}_{\frac{\infty}{2}}$-modules  $L_w\otimes L_w$ with the
highest weight vectors $w\cdot v$ $({\overline w}=1)$ are irreducible,
and it follows from Ch.~1 (relation~(7))
that these vectors are not $\widehat{\g}_{\infty ,0}(\C \
[t]/t^2)$-singular in $Ind_{\mu_1,\mu_t,s}$. On the other hand,
it follows that the intersection of the greatest proper submodule in
$Ind_{\mu_1,\mu_t,s}$ with $S^*(\A_-)\otimes 1\subset S^*(\A_-)\otimes
S^*(\overline{\A_-})$ is zero. This proves Lemma~7 since
$S^*(\A_-)=\bigoplus\limits_wL_w\otimes L_w$ as a
$\g^{(1)}_{\frac{\infty}{2}}
\oplus \g^{(2)}_{\frac{\infty}{2}}$-module, see Ch.~1, Sec.~1.
\enddemo
\remark
{4.1.4. Remarks}

1. The vectors ${\overline w}\cdot v$ are obviously singular vectors
of the $\widehat{\g}_{\infty ,0}(\C \ [t]/t^2)$-module
$Ind_{\mu_1,\mu_t,s}$ and therefore, the vectors
$w\cdot {\overline w}\cdot v$ $({\overline w}\ne 1)$ and the irreducible
$\g^{(1)}_{\frac{\infty}{2}}\oplus \g^{(2)}_ {\frac{\infty}{2}}$-modules
generated by them belong to the maximal submodule in
$Ind_{\mu_1,\mu_t,s}$.

2. Clearly, the $\widehat{\g}_{\infty ,0}(\C \ [t]/t^2)$-module
$Ind_{\mu_1,\mu_t,s}$ does not depend on $\mu_1$ in the obvious sense.
\endremark
\subsubhead
4.2
\endsubsubhead
Recall that $\g(\lambda )={\Lie}\ ({\Cal U}_{\lambda})$,
where ${\Cal U}_{\lambda}=U(\s_2)\bigg/ \biggl( \Delta -
{\displaystyle\frac{\lambda (\lambda +2)}{2}}\biggr)$, and $\Delta$
is the Casimir operator
$e\cdot f+f\cdot e+{\displaystyle\frac{h\cdot h}{2}}\in
U(\s_2)$, $\lambda \in \C$. Being a vector space,
$$
\align
{\Cal U}_{\lambda}=& \C \ [h]\oplus \\
& \oplus e\cdot \C \ [h]\oplus
e^2\cdot \C \ [h]\oplus \ldots \\
& \oplus f\cdot \C \ [h]\oplus
f^2\cdot \C \ [h]\oplus \ldots
\endalign
$$
($\{ e,f,h\}$ is a standard basis in  $\s_2$). The relations
$$
\left\{
\aligned
h\cdot e&=e\cdot (h+2)\\
h\cdot f&=f\cdot (h-2)\\
e\cdot f&=T_1(h)={\displaystyle\frac{1}{2}}\biggl(
h-{\displaystyle\frac{h^2}{2}}+{\displaystyle\frac{\lambda (\lambda
+2)}{2}}\biggr)\\
f\cdot e&=T_2(h)={\displaystyle\frac{1}{2}}\biggl(
-h-{\displaystyle\frac{h^2}{2}}+{\displaystyle\frac{\lambda (\lambda
+2)}{2}}\biggr) .
\endaligned
\right.
\tag 32
$$
\subsubhead
4.2.1
\endsubsubhead
We define the $\M^s_{\infty}(\C \ [t]/t^m)$ algebra by relations~(30)
(setting $A=\prod\limits_{i\in
\Z}(\C \ [t]/t^m)_i$) and replacing the last two conditions in~(30) by
$$
\left\{
\aligned
&.\ .\ .\ .\ .\ .\ .\ .\ .\ .\ .\ .\ .\ .\ .\ .\ .\ .\ .\ .\ .\ .\ .\ .\ .\\
&e\cdot f=\biggl( T_1(s-2i+t),\ i\in \Z \biggr) \in A\\
&f\cdot e=\biggl(
T_2(s-2i+t),\ i\in \Z \biggr) \in A,
\endaligned
\right.
\tag 33
$$
$$
m\in \Z_{\geq 1},\ s\in \C \ .
$$

In order to construct the mapping ${\Cal U}_{\lambda}$ in
$\M^s_{\infty}(\C [t]/t^m)$, it suffices to determine the mapping
$\varphi_s: \C \ [h]\to A$
such that the diagram
$$
\CD
{\C}\ [h] @>{\varphi_s}>>  A\\
@V\sigma VV @VVDV\\
{\C}\ [h] @>{\varphi_s}>>  A
\endCD
\tag 34
$$
is commutative $\bigl( \sigma (p(h))=p(h+2)\bigr)$. Then $\varphi_s$
is extended to ${\Cal U}_{\lambda}$, being the homomorphism of an
associative algebras, with the use of the relations
$\varphi_s(e)=e,\ \varphi_s(f)=f$.

It immediately follows from~(34) that $\varphi_s(p(h))=\bigl( p(s-2i+t),\
i\in \Z \bigr) \in A$. Obviously, $\varphi_s: {\Cal U}_{\lambda}\mapsto
\M^s_{\infty}(\C \ [t]/t^m)$ is an embedding.
\subsubhead
4.2.2
\endsubsubhead
Let $m=2$.

\proclaim{Lemma~8}
For a general $\lambda$ we have
\rom {(i)} $\M^s_{\infty}(\C \ [t]/t^2)\cong \M_{\infty}(\C \ [t]/t^2)$
if $T_1(s-2i)\ne 0$ for all $i\in \Z$,

\rom {(ii)} $\M^s_{\infty}\cong \M_{\infty ,i}$ if $T_1(s-2i)=0$
\endproclaim
\demo{Proof} This is obvious; see also Lemma~3.
\enddemo

We denote by $\g^s_{\infty}(\C \ [t]/t^m)$ the Lie algebra constructed with
the use of the associative algebra\newline
$\M^s_{\infty}(\C \ [t]/t^m)$.

The cocycle $\alpha$ on the Lie algebra
$\g^s_{\infty}(\C \ [t]/t^2)$ with the values in $\C \ [t]/t^2$
is constructed in the standard way,
$$
\aligned
\alpha (E_{ij},E_{ji})&={\cases
P_{j-i}(s-2i)& \text{for}\ i\leq 0,\ j\geq 1\\
0,& \text{otherwise}
\endcases}\\
\alpha (E_{ij},E_{kl})&=0\ \text{for}\ i\ne l\ \text{and}\ j\ne k
\endaligned
\tag 35
$$
and is continued by $\C \ [t]/t^2$-linearity (here
$P_l(h)=e^l\cdot f^l$, $E_{ij}=1_i\cdot e^{j-i}$ for $i<j$.) In
case~(i) of Lemma~8, this cocycle is the inverse image of the
standard cocycle on $\g_{\infty}(\C \ [t]/t^2)$ and, in case~(ii)
of Lemma~8, the corresponding central extension is the direct sum of the
trivial and the one-dimensional extension.
\subsubhead
4.2.3
\endsubsubhead
We set $h_k=[e,fh^{k-1}]$ in $\g(\lambda )$ and then, in the base
$\{ h_k\}$, the highest weight of the representation induced from the
parabolic subalgebra $\p_{\alpha}$ is equal to
$$
\left\{
\aligned
&\chi (h_1) = \chi (h)\\
&\chi (h_2) = \alpha \chi (h)\\
&\chi (h_3) = \alpha^2 \chi (h)\\
&.\ .\ .\ .\ .\ .\ .\ .\ .\ .\ .
\endaligned
\right.
\tag 36
$$
Next, the highest weight of the representation
$\theta^*_s(Ind_{\mu_1,\mu_t,s})$
($\theta_s: \g(\lambda )\hookrightarrow \widehat{\g}^s_{\infty}(\C \
[t]/t^2)$) is equal, by relation~(35), to
$$
\left\{
\aligned
&\chi (h_1) = T_1(s)\cdot \mu_1+T'_1(s)\mu_t\\
&\chi (h_2) = T_1(s)\cdot (s-2)\mu_1+\bigl( T'_1(s)\cdot
(s-2)+T_1(s)\bigr) \mu_t\\
&\chi (h_3) = T_1(s)\cdot (s-2)^2\mu_1+\bigl( T'_1(s)\cdot
(s-2)^2+2(s-2)\cdot T_1(s)\bigr) \mu_t\\
&\chi (h_4) = T_1(s)\cdot (s-2)^3\mu_1+\bigl( T'_1(s)\cdot
(s-2)^3+3(s-2)^2\cdot T_1(s)\bigr) \mu_t\\
&.\ .\ .\ .\ .\ .\ .\ .\ .\ .\ .\ .\ .\ .\ .\ .\ .\ .\ .\ .\ .\ .\ .\ .\ .\
.\ .\ .\ .\ .\ .\ .\ .\ .\ .\ .\ .\ .\ .\ .\ .\ .\ .\ .\ .\ .\ .
\endaligned
\right.
\tag 37
$$

Under the embedding $\theta_s:\g(\lambda )\hookrightarrow
\g^s_{\infty}(\C \ )$ the highest weight of the module $\theta^*_s(\I)$
is equal to
$$
\left\{
\aligned
&\chi (h_1) = T_1(s)\cdot \mu \\
&\chi (h_2) = (s-2)\cdot T_1(s)\cdot \mu \\
&.\ .\ .\ .\ .\ .\ .\ .\ .\ .\ .\ .\ .\ .\ .\ .\ .\ .\ .
\endaligned
\right.
\tag 38
$$
Comparing~(36) with (38), we see that if $T_1(\alpha +2)=0$
the corresponding highest weight exists in~(36) and does not exist in~(38)
for all $\chi (h)$ (since  $s=\alpha +2$).

We set $T_1(s)=0$ in (37);
then we can easily pass from~(37) to (36) by setting
$\alpha =s-2$, $T_1(\alpha +2)=0$,
$\mu_1$ being arbitrary, $\chi (h)=T'_1(s)\cdot \mu_t$. {\it The main fact
here is that \rom(for the general $\lambda$\rom)
$T_1(h)$ has no multiple roots, i.e.,} $T_1(s)=0\Rightarrow
T'_1(s)\ne 0$.
\subsubhead
4.2.4
\endsubsubhead
Thus, $T_1(s)=0$, and we are in the situation of~(ii) of Lemma~8,
$\widehat{\g^s_{\infty}}(\C \ [t]/t^2)\cong
\widehat{\g}_{\infty ,0}(\C \ [t]/t^2)$. Lemma~7 states that the character
of the irreducible quotient module of the $\widehat{\g}_{\infty ,0}(\C \
[t]/t^2)$-module  $Ind_{\mu_1,\mu_t,s}$ for $\mu_t\ne 0$ is greater
or equal to $\prod\limits_{i\geq 1}(1-q^i)^{-i}$ and we need the following
lemma.

\proclaim{Lemma~9~\rm[2]}
Any $\g(\lambda )$-invariant subspace of the module
$Ind_{\mu_1,\mu_t,s}$ is invariant with respect to
$\widehat{\g^s_{\infty}}(\C \ [t]/t^2)$ (for all $s\in
\C$).
\endproclaim

\demo{Proof}
See the corollary  of Lemma~6.
\enddemo

Lemma~7, Lemma~8, and Lemma~9 imply Theorem~3 (see Subsec.~4.0).
\remark
{4.2.5. Remarks}

1. It follows from relation~(37) that all representations of
$gl(\lambda )$ induced from the parabolic subalgebras
$\p_{\alpha ,\alpha}$, $\alpha\in \C$ (generated by
$\n_+$, $\h$,  and $\{ f(h-\alpha )^2p(h),\ p(h)\in
\C \ [h]\}$) have the form
$\theta^*_s(Ind_{\mu_1,\mu_t,s})$.

Indeed, we have the recurrent condition
$$
\chi (h_k)=2\alpha \chi (h_{k-1})-\alpha^2\chi (h_{k-2});
$$
imposed on the highest weight $\chi$ of the representation induced from
$\p_{\alpha ,\alpha}$,  $\alpha^k$, and
$k\alpha^{k-1}$ are the corresponding basic solutions.

2. Let $T_1(s)=0$ and $T_1(h)=a(h-s)+b(h-s)^2\ (a\ne 0)$.

We set
$$
\mu_1=\frac{1}{a}\cdot \frac{1}{(h-s)^2};\quad \mu_t=-\frac{1}{b}\cdot
\frac{1}{(h-s)}.
$$
Then, as $h\to s$, relations~(37) give the {\it limiting\/}
values, and we find that for $T_1(\alpha +2)=0$, the representations
induced from $\p_{\alpha ,\alpha}$ do not have notation~(37).
We could have obtained them, by analogy with Subsec.~4.2.3, as the inverse
images under the inclusion of $\g(\lambda )$ into
$\widehat{\g^s_{\infty}}(\C \ [t]/t^3)$.

Next, for obtaining identity~(1) we used in Sec.~3 not the whole
expression for the determinant of Shapovalov's form on the levels
of the $\widehat{\g}_{\infty}$-module of $\I$ but only its {\it degrees\/}
with respect to $\mu$ on the levels. However, the ``highest'' local
identities obtained in this way are simply equal to the degrees of
identity~(1) (see also Appendix~B).
\endremark

\head
Chapter 3.\\
The Lie Algebra of Functions on a Hyperboloid:\\
the Global Identity
\endhead

\def\BI{\Bbb I}
\def\D{\operatorname{Det}}
\def\n{{\frak n}}
\def\h{{\frak h}}
\def\p{{\frak p}}
\def\G{{\frak g}}
\def\A{{\frak A}}
\def\I{Ind_{\mu}}
\def\In{Ind_a}
\def\g{gl}
\def\e{\operatorname{deg}}
\def\s{sl}

\subhead
1. The Flat Deformation in $\g(\lambda )$
and the Principal Vector Bundle  on $S^2$
\endsubhead
\subsubhead
1.1
\endsubsubhead
If we assume that the Lie algebra $\g(\lambda )$ corresponds to the
point $\lambda \in \C =S^2\setminus \{
\infty \}$, then the Lie algebra of functions on a hyperboloid corresponds
to the point $\{ \infty \}\in S^2$. To be more precise, we shall consider the
standard symplectic foliation over $\s_2(\C)^*$. We call the general symplectic
leaf of this foliation, which is a nondegenerate quadratic surface in $\C^3$,
a hyperboloid. The induced Poisson bracket defines the structure of the Lie
algebra on regular functions on a symplectic leaf. These Lie algebras are
isomorphic for all nondegenerate symplectic leaves.

There exists a flat deformation of the Lie algebras of functions
on a hyperboloid  into the Lie algebra $\g(\lambda_t)$ in the
neighborhood of $t=0$, where
$\lambda_t\cdot (\lambda_t+2)={\displaystyle\frac{\lambda
(\lambda +2)}{t^2}}$ and $\lambda \in \C \setminus \{
0\}$ is fixed.  We fix $t\ne 0$. We shall consider the Lie algebra
$\s_2(\C)$ as an algebra with the generators $e$, $h$, $f$, and
relations $[e,f]=th$, $[h,e]=2te$, $[h,f]=-2tf$. We denote
by $U_{\lambda}(t)$ the quotient algebra of the corresponding universal
enveloping algebra with respect to the ideal generated by the
{\it central\/} element
$$
2ef - th + \frac{h^2}{2} - \frac{\lambda (\lambda +2)}{2},\quad
\text{for}\ t\ne 0,
$$
$$
U_{\lambda}(t)\cong U_{\lambda_t}=U(\s_2)\bigg/ \biggl( \Delta -
\frac{\lambda (\lambda +2)}{t^2}\biggr) ,
$$
where $\Delta$ is the standard Casimir operator in $U(\s_2)$. The corresponding
Lie algebra with the bracket $[a,b]=a\cdot b-b\cdot a$ is isomorphic
to the Lie algebra
$\g(\lambda_t)$. We denote by ${\Lie}'U_{\lambda}(t)$ the corresponding
Lie algebra with the bracket
$[a,b]={\displaystyle\frac{a\cdot b-b\cdot a}{t}}$
The Lie algebra ${\frak s}=\lim\limits_{t\to 0}{\Lie}'U_{\lambda}(t)$
is isomorphic to the Lie algebra of regular functions on the symplectic
leaf $\biggl\{
2ef+{\displaystyle\frac{h^2}{2}}={\displaystyle\frac{\lambda (\lambda
+2)}{2}}\biggr\}$ of the foliation on $\s^*_2$ with an induced Poisson bracket.
(If we simplify this construction without factorizing with respect to the
central element, then we get a Poisson bracket on regular
functions on $\s^*_2$, i.e., on $\C \ [e,f,h]$.)
\subsubhead
1.2
\endsubsubhead
We denote $U=S^2\setminus \{ \infty \}$ and $V=S^2\setminus \{
0\}$ and assume that the parameter $t$ runs over all points from
$S^2$. Thus we described in Subsec.~1.1 what we have in $U$.
We denote by $e_U$, $f_U$, $h_U$ the elements $e$, $f$, $h$ from Subsec.~1.1
respectively.
Suppose now that the Lie algebra
$\g(\lambda_t)$ corresponds to the point $t\in V=S^2\setminus \{ 0\}$
($\g(0)$ corresponds to the point $\{ \infty \}$).
We denote the elements $e$, $f$, $h$ appearing in the definition of
$\g(\lambda_t)$ by $e_V$, $f_V$, $h_V$ respectively.
We want to calculate the corresponding transition functions. First, we note
that the Lie subalgebra $\s_2$ is globally defined, i.e. the transition
functions are identical. In general, the isomorphism
$\alpha : U_{\lambda}(t)\to U_{\lambda}$
is defined as $\alpha : e^kh^l\mapsto t^{k+l}\cdot e^kh^l$ etc, and
$\beta : \xi \mapsto {\displaystyle\frac{\xi}{t}}$ under the isomorphism
$\beta : {\Lie}'U_{\lambda}(t)\to {\Lie}U_{\lambda}(t)$
for any $\xi$. Therefore on $U\cap V$ we have
$$
\frac{e^k_Uh^l_U}{e^k_Vh^l_V} = t^{k+l-1}.
$$
Thus the family of our Lie algebras on $S^2$ is decomposed into the
(infinite) direct sum of line bundles on $S^2$.
\subsubhead
1.3
\endsubsubhead
At every point $t\in S^2$ we choose a parabolic subalgebra of the
corresponding Lie algebra which contains $\n_+$ holomorphically
dependent on $t$.
We shall only consider parabolic subalgebras of the
``maximal size'', i.e. subalgebras such that the first level of the
induced representation is one-dimensional. Recall that for a fixed $t$
these subalgebras depend on one complex parameter $\alpha$ and the
subalgebra ${\frak p}_{\alpha}$ is generated by $\n_+$, $\h$, and
$$
\n_-^{\alpha}=\biggl\{ p(h)\cdot (h-\alpha )\cdot f,\ p(h)\in
\C \ [h]\biggr\}.
$$
Furthermore, for every $t\in S^2$ we consider the representation
induced from the corresponding parabolic subalgebra for a certain
value of the highest weight $\chi_t(h)$ for all $t\in S^2$. Then we can
regard the spaces of representations as the direct sum of line bundles
on $S^2$  and
$$
\sum\limits_{k\geq 0}\left( \
\sum\limits^{n_k}_{l=1}a^{C_{l,k}}\right) \cdot q^k=
{\displaystyle\frac{1}{(1-q)}} \cdot
{\displaystyle\frac{1}{(1-aq^2)(1-a^2q^2)}} \cdot
{\displaystyle\frac{1}{(1-a^2q^3)(1-a^3q^3)(1-a^4q^3)}} \ldots \ .
\tag 1
$$
Here $n_k$ is the dimension of the $k$th level of the induced representation
and $C_{1,k};\ldots ;C_{n_k,k}$ are Chern classes of the line bundles into
whose direct sum the $k$th level of the representation is decomposed.
(Every monomial $f^ph^q$ has a Chern class equal to
$p+q-1$ in addition to its
grading equal to $p$. In this sence, the last relation is the refinement of our
relation for the character $\prod\limits_{i\geq
1}{\displaystyle\frac{1}{(1-q^i)^i}}$ with due account of the Chern class.)

The quadratic form, namely, Shapovalov's form, is invariantly defined
in every fiber of a definite level. Clearly, its determinant is not
uniquely defined. However, we choose a holomorphic frame of the bundle of the
$k$th level subject to the decomposition of this bundle into a direct sum.
Suppose that $U$, $V$ are our coordinate maps on
$S^2$ and $\xi_{1,U};\ldots ;\xi_{n_k,U}$
and $\xi_{1,V};\ldots ;\xi_{n_k,V}$ are the corresponding sections in these
coordinates. Then
${\displaystyle\frac{\xi_{i,U}}{\xi_{i,V}}}=t^{C_{i,k}}$ and
$$
\frac{\det \lnorm \langle \xi_{i,U}, \xi_{j,U}\rangle ,\
i,j=1,\ldots , n_k\rnorm}{\det \lnorm \langle \xi_{i,V}, \xi_{j,V}\rangle
,\ i,j=1,\ldots , n_k\rnorm} = t^{2\sum\limits^{n_k}_{i=1}C_{i,k}}.
$$
Therefore, the determinant of Shapovalov's form is a
({\it not uniquely\/} defined) holomorphic section of the line bundle over
$S^2$ with the Chern class $2\sum\limits^{n_k}_{i=1}C_{i,k}$
(on the $k$th level).
We denote this Chern class by $C_k$. Then
$$
\sum\limits_{k\geq 1}C_k\cdot q^k =
{\displaystyle\frac{d}{da}}\Biggl[ {\displaystyle\frac{1}{(1-q)}}\ \cdot
{\displaystyle\frac{1}{(1-a^2q^2)(1-a^4q^2)}}\ \cdot \
{\displaystyle\frac{1}{(1-a^4q^3)(1-a^6q^3)(1-a^8q^3)}}\ \ldots \Biggr]
\Bigg\vert_{a=1}.
\tag 2
$$
This is a simple consequence of relation~(1).
\subsubhead
1.4
\endsubsubhead
We choose the subalgebras ${\frak p}_{\alpha ,t}$ and characters
$\chi_t(h)$ for all $t\in S^2$ such that for $t=0$ (corresponding
to the Lie algebra of functions on a hyperboloid)
$\frak p$ and $\chi (h)$ are in the
{\it general position}. In this case we can calculate the Chern class of
the bundle of the determinant of Shapovalov's form in two ways, namely,
as in Subsec.~1.3 and as the sum of zeros of the determinant of
Shapovalov's form with multiplicities (the multiplicities of zeros are
uniquely defined).
In the next section, we shall calculate the Chern class with the second
technique, using the results of Chs.~1 and 2 on the assumption of the
irreducibility  for $t=0$.
Equating the corresponding generating function to expression~(2),
we obtain an identity with power series which is equivalent to
the irreducibility for the general parameters of the induced
representations of the Lie algebra of functions on a hyperboloid.
We shall prove a more exact result in Sec.~3 (see Theorem~1).
\subhead
2. The Main Calculation
\endsubhead
\subsubhead
2.1
\endsubsubhead
\proclaim{Lemma~1}
Let $\alpha$, $\chi \in \C$ be arbitrary. We set
$\alpha_t={\displaystyle\frac{\alpha}{t}},\ \chi_t(h_V)=\chi$ in the
coordinates on $V=S^2\setminus \{ 0\}$.
Then the corresponding limiting parameters on
$U=S^2\setminus \{\infty \}$ as $t\to 0$, are
equal to $\alpha_0=\alpha$ and $\chi_0(h_U)=\chi$
\rom(in the coordinates on $U$\rom{).}
\endproclaim
\demo{Proof}
The direct verification.
\enddemo

Thus, for any parameters of the induced representation for
$t=0$, we have a holomorphic continuation to the whole sphere.

{\it Everywhere in what follows  $\alpha$ and $\chi$ are considered to be
of the general position}.

We shall need relations that connect the $(s,\mu )$-parametrization
of the induced representations of $\g(\lambda )$ with the one-dimensional
first level with the $(\alpha ,\chi (h))$-parametrization (see Subsec.~2.4
of Ch.~2).  We have
$$
\left\{
\aligned
&T_1(s)\mu =\chi (h)\\
&(2s-2)T_1(s)\mu =2(1+\alpha )\chi (h).
\endaligned
\right.
\tag 3
$$
The meaning of this is that in the $(s,\mu )$-parametrization the
reducibility condition can be naturally formulated:
$$
\left.
\alignedat 2
&\text{(i)}&\quad &\mu \in \Z \\
&\text{(ii)}&\quad& T_1(s+2i)=0\ \text{for}\\
&&\quad& \text{a certain}\ i\in \Z.
\endalignedat
\right\}
\tag 4
$$
(recall that
$T_1(s)={\displaystyle\frac{1}{2}}\Bigl(
s-{\displaystyle\frac{s^2}{2}}+{\displaystyle\frac{\lambda
(\lambda +2)}{2}}\Bigr)$. For $T_1(s)\ne 0$ it follows from~(3) that
$$
s = \alpha +2,\ \text{and}\ \mu = \chi (h)\big/ T_1(\alpha +2).
$$

Let $\chi (h)$ be arbitrary and  $T_1(\alpha +2)=0$. In this case,
the corresponding representation does not have the
$(s,\mu )$-parametrization. Thus, under the conditions of Lemma~1, we have
on $V=S^2\setminus \{ 0\}$
$$
\left.
\alignedat 2
&\text{(i)}&\quad \chi (h)\bigg/ T_{1,t}&\biggl(
{\displaystyle\frac{\alpha}{t}}+2\biggr) \in \Z ;\\
&\text{(ii)}&\quad
T_{1,t}\biggl( {\displaystyle\frac{\alpha}{t}}+2i\biggr)& =0\
\text{for a certain}\\
&&\quad& i\in \Z ;\\
\endalignedat
\right\}
\tag 5
$$
$$
T_{1,t}(s) ={\displaystyle\frac{1}{2}}\biggl(
s-{\displaystyle\frac{s^2}{2}}+{\displaystyle\frac{\lambda
(\lambda +2)}{2t^2}}\biggr).
$$
\subsubhead
2.2
\endsubsubhead
Thus, all points $t\in S^2\setminus \{ 0\}$ at which the determinant
of Shapovalov's form has a zero are given by formulas~(5), (i) and (ii).
In addition, it follows from Theorem~1, Sec.~3 (see below) that for $t=0$
the determinant of Shapovalov's form does not vanish (for the general
$\alpha,\chi (h)$).

Our main purpose is to determine the {\it multiplicities\/}
of zeros according as $t\in S^2$. For example, in case~(5), (i)
we must connect the {\it multiplicity of the zero of the determinant\/}
of Shapovalov's form for
$\widehat{\g}_{\infty}$ for the integral central charge $\mu$,
{\it as a function of $\mu$, found in Ch.~1 with the multiplicity of the zero
as a function of $t\in S^2$}. In fact, the results of Ch.~1 and 2
are sufficient for all cases~((5), (i) and (ii)).

Recall (see~(29) Ch.~2) that the determinant of Shapovalov's form
of the $\g(\lambda )$-module $\theta^*_s(Ind_{\mu ,s})$ as a function of
$s$ and $\mu$, is given by the relation
$$
\det \Phi_{\g(\lambda )}(s,\mu ) = \Biggl(
\prod\Sb\text{over all}\\
w\text{from}\\
\text{the $k$th level}\endSb
\prod\limits^p_{m=1}\prod\limits^{i_p-1}_{l=0}
T_1(s-2j_m-2l)\Biggr) \times
\det \Phi_{\widehat{\g}_{\infty}}(\mu ).
$$

Here we set
$$
w=\bigl( f^{i_1}\delta_{j_1,S}\bigr) \cdot \ldots \cdot \bigl(
f^{i_p}\delta_{j_p,S}\bigr) \cdot v,
$$
where $i_1+\ldots +i_p=k$ and
$$
j_q\in [-i_q+1,0].
\tag 7
$$
\subsubhead
2.3
\endsubsubhead
Here we find the multiplicities of zeros $t\in S^2\setminus \{ 0\}$
in cases~(5), (i) and (ii).
\subsubhead
2.3.1
\endsubsubhead
In case~(5), (i) on the assumption that $\alpha$ and $\lambda$ are general
$T_{1,t_0}\biggl( {\displaystyle\frac{\alpha}{t_0}}+2i\biggr)\ne 0$
is nonzero for any $i\in \Z$, and therefore the product in~(6)
does not vanish. Thus, if $\mu_0\in \Z$, then the multiplicity of the root
$t_0$ according as $t\in S^2$ in the determinant of Shapovalov's form of the
$k$th level for
$\g(\lambda_t)$ is equal to the multiplicity of the root $\mu_0$ according as
$\mu \in \C$ in the determinant of Shapovalov's form of the $k$th level of
$\widehat{\g}_{\infty}$.
\subsubhead
2.3.2
\endsubsubhead
Case~(5), (ii) for $i=1$ was considered in Ch.~2, namely, in this case the
corresponding representation $(\alpha ,\chi (h))$ does not have the
$(s,\mu )$-parametrization, and
$t=t_0\ \mu \sim {\displaystyle\frac{1}{t-t_0}}$ in the deleted neighborhoods.
Therefore, here, in relation~(6),
$$
\det \Phi_{\widehat{\g}_{\infty}}(\mu )\sim \frac{1}{(t-t_0)^{\deg_k}},
$$
where $\deg_k$ is the degree of  $\det \Phi_{\widehat{\g}_{\infty}}(\mu )$ with
respect to $\mu$ on the $k$th level.
On the other hand, according to condition~(7), the multiplicity of zero in~(6)
for the general $\alpha$, $\lambda$ is equal to $p$ for every
$$
w=\bigl( f^{i_1}\delta_{j_1,S}\bigr) \cdot \ldots \cdot \bigl(
f^{i_p}\delta_{j_p,S}\bigr) \cdot v,
$$
since $T_1(s-2i)$ is zero only for $i=0$.

Thus, in accordance with Lemma~12 from Subsec.~2.2 of Ch.~1, we see that the
irreducibility of the representation of the Lie algebra $\g(\lambda )$,
induced from ${\frak p}_{\alpha}$, for $T_1(\alpha +2)=0$ and the general
$\lambda$, is equivalent to the following local identity from Ch.~2:
$$
\frac{d}{da}\Biggl( \ \prod_{i\geq 1}\frac{1}{(1-aq^i)^i}\Biggr)
\Bigg\vert_{a=1} = \sum\Sb\text{over all}\\
w=\D_1^{l_1}\cdot \ldots \cdot
\D_k^{l_k}\endSb\# \ D(w)\cdot q^{\sum l_i\cdot i^2}\cdot \biggl(
\chi (w)\biggr)^2,
\tag 8
$$
where $\# D(w)$ is the number of squares in the Young diagram corresponding
to $w$  (see Introduction, Fig.~1) and $\chi (w)$ is the corresponding semiinfinite character
(see Subsec.~0.11 or Ch.~1).
\subsubhead
2.3.3
\endsubsubhead
Let us consider case~(5), (ii) for $i\ne 1$. In this case, the corresponding
representation (for the general $\alpha$, $\lambda$) belongs to the
$(s,\mu )$-parametrization and the part of relation~(6) dependent on
$\mu$ does not have any zeros of poles. Furthermore, let
$\lmod i-1\rmod =k_+$. Then the generating function
$$
\sum_{k\geq 0}\left.\cases
\text{the multiplicity of zero in the neighborhood of}\ t=t_0\\
\text{of the determinant of Shapovalov's form on the}\\
k\text{th}\ \text{level}\endcases\right\}\cdot q^k
$$
is equal, by virtue of relations~(6), (7), to
$$
\frac{d}{da}\Biggl[ \ \prod^{k_+}_{i=1}\frac{1}{(1-q^i)^i}\cdot
\prod^{\infty}_{i=k_++1}\frac{1}{(1-q^i)^{k_+}\cdot
(1-aq^i)^{i-k_+}}\Biggr] \Bigg\vert_{a=1}.
\tag 9
$$
\subsubhead
2.4
\endsubsubhead
Here we shall obtain a global identity equivalent to the theorem on the
irreducibility, for the general parameters, of the induced representation
of the Lie algebra of functions on a hyperboloid
(with the one-dimensional first level).

None of the values of $t\in S^2\setminus \{ 0\}$ meets condition~(5),
(i) for $\mu =0$ and a one pair  of points on $S^2\setminus \{ 0\}$ meets it
for every $\mu \in \Z \setminus \{ 0\}$.

According to Subsec.~2.3.1, the corresponding generating function for the
sum of multiplicities of zeros of the determinant of Shapovaslov's form
of the $k$th level for all $t\in S^2\setminus
\{ 0\}$ that satisfy condition~(5), (i) with $\mu \ne 0$ is equal to
$$
2\sum_{a\in \Z \setminus \{ 0\}}\left.\cases
\text{the multiplicity of the root}\ \mu =a\\
\text{in the determinant}\\
\text{of Shapovalov's form for}\ \widehat{\g}_{\infty}\ \text{on the}\\
k\text{th\  level}
\endcases\right\} \cdot q^k.
\tag 10
$$

For every $k_+\ne 0$ (see Subsec.~2.3.3) condition~(5), (ii) is satisfied
by {\it four} points on  $S^2\setminus \{ 0\}$ and for
$k_+=0$ it is satisfied by {\it two} points.

It follows from Theorem~1, Sec.~3 (see belolw) that for the general
$\alpha$, $\lambda$ the corresponding induced representation of the Lie
algebra of functions on a hyperboloid is irreducible. Then the Chern
class found in~(2) is equal to the sum of multiplicities of zeros of the
determinant of Shapovalov's form {\it for\/} $t\in S^2\setminus \{ 0\}$.
We have
$$
\gather
\frac{d}{da}\biggl[
\frac{1}{(1-q)}\ \cdot \ \frac{1}{(1-a^2q^2)(1-a^4q^2)}\ \cdot
\ \frac{1}{(1-a^4q^3)\ldots}\ \ldots \biggr] \bigg\vert_{a=1}=\\
1+2\sum_{k\geq 1}\sum_{a\in \Z \setminus \{
0\}}\left.\cases \text{the multiplicity of the root}\ \mu =a\ \text{
in the determinant}\\ \text{of Shapovalov's form of}\ \widehat{\g}_{\infty}\
\text{on the}\ k\text{th level}
\endcases\right\} \cdot q^k +\\
2\frac{d}{da}\Biggl[ \ \prod_{i\geq 1}\frac{1}{(1-aq^i)^i}\Biggr]
\Bigg\vert_{a=1}-2\sum_{k\geq 1}\left.\cases
\text{the degree with respect to}\ \mu \ \text{of the determinant}\\
\text{of Shapovalov's form of}\
\widehat{\g}_{\infty}\ \text{on the}\ k\text{th}\\
\text{level}
\endcases\right\} \cdot q^k +\\
4\sum_{k_+\geq 1}\frac{d}{da}\Biggl[ \
\prod^{k_+}_{i=1}\frac{1}{(1-q^i)^i}\cdot
\prod^{\infty}_{i=k_++1}\frac{1}{(1-q^i)^{k_+}(1-aq^i)^{i-k_+}}\Biggr]
\Bigg\vert_{a=1}.
\tag 11
\endgather
$$
In accordance with the local identity from Subsec.~2.3.1, the sum of the
two last but one terms is equal to 0. However, in form~(11) we easily
arrive at the following form of the global identity.

{\it Global identity}:
$$
\gather
\frac{d}{da}\biggl[ \frac{1}{(1-q)}\ \cdot \ \frac{1}{(1-aq^2)(1-a^2q^2)}\
\cdot \ \frac{1}{(1-a^2q^3)(1-a^3q^3)(1-a^4q^3)}\ \ldots \biggr]
\bigg\vert_{a=1}=\\
\frac{d}{da}\Biggl[ \ \prod_{i\geq 1}\frac{1}{(1-aq^i)^i}\Biggr]
\Bigg\vert_{a=1}+\\
2\sum_{k_+\geq 1}\frac{d}{da}\Biggl[ \
\prod^{k_+}_{i=1}\frac{1}{(1-q^i)^i}\cdot
\prod^{\infty}_{i=k_++1}\frac{1}{(1-q^i)^{k_+}\cdot
(1-aq^i)^{i-k_+}}\Biggr] \Bigg\vert_{a=1}-\\
\sum\Sb\text{over all}\\
\text{Young diagrams}\\
w=\D_1^{l_1}\cdot \ldots \cdot \D_k^{l_k}\endSb
\left.\cases
\text{the multiplicity of the root}\ \mu =0\\
\text{in}\ p_{w}(\mu )=\langle w,w\rangle
\endcases\right\} \cdot
q^{\sum l_i\cdot i^2}\cdot \biggl( \chi (w)\biggr)^2.
\tag 12
\endgather
$$

According to relation~(8) from Ch.~1, the multiplicity of the root
$\mu =0$ in $p_w(\mu)=\langle w,w\rangle$ is equal to the length of the
``central'' diagonal in the corresponding Young diagram that starts from
the upper left vertex (see Fig.~1, Introduction).

Recall the definition of the semi infinite character $\chi (w)$.
Let us consider the Lie algebra of finite matrices
$(a_{ij}),\ i,j=1\ldots \infty$. Suppose that
$\alpha_1^{\vee},\alpha_2^{\vee},\ldots$ are the corresponding simple coroots.
Then $\chi (w)$ is the character of the irreducible representation of this
Lie algebra with the highest weight
$\chi: \chi (\alpha_1^{\vee})=l_1,\ldots ,\chi (\alpha_k^{\vee})=l_k,\ \chi
(\alpha_{k+1}^{\vee})=\ldots =0$. Thus
$\chi (w)={\displaystyle\frac{1}{(1-q)\cdot \ldots \cdot
(1-q^k)}}$ for $w=\D_k$ as well as for
$w=\D_1^k$,
\subhead
3. Proof of the Global Identity: Reduction to the Nilpotent Case
\endsubhead
\subsubhead
3.0
\endsubsubhead
According to the results of Sec.~2, the global identity~(12)
is equivalent to the irreducibility of representations
of the Lie algebra $\frak s$ of functions on the hyperboloid
$\bigl\{2ef+{\displaystyle\frac{h^2}{2}}={\displaystyle\frac{\lambda
(\lambda +2)}{2}}\bigr\}$, $\lambda (\lambda +2)\ne 0$ induced from
the parabolic subalgebras $\p_{\alpha}$ with the highest weight
$\chi$ for the general values of $\alpha$ and $\chi(h)$. We denote by
$\frak s_0$  the Lie algebra of functions on the cone
$\bigl\{ 2ef+{\displaystyle\frac{h^2}{2}}=0\bigr\}$ which is a degenerate
symplectic leaf of the foliation in $\s_2(\C \ )^*$. It is clear that
it suffices to prove a similar statement for
${\frak s}_0$ (of the Lie algebra of functions on all nondegenerate symplectic
leaves are isomorphic). But the Lie algebra
${\frak s}_0$ is nilpotent with an accuracy to within the subalgebra
$\s_2$ contained in it, and we can apply Kirillov's theory to it
(see Appeendix~A).
\subsubhead
3.1
\endsubsubhead
The following theorem is the main result of this section.

\proclaim{Theorem~1}
The representation of the Lie algebra ${\frak s}$ induced from the
parabolic subalgebra $\p_{\alpha}$ with the highest weight
$\chi \in \h^*$ is reducible for

\rom{(1)} $\chi (h) =0$\newline
or

\rom{(2)} $\alpha^2=\lambda (\lambda +2)$
and ireducible otherwise.
\endproclaim
\remark
{Remarks}

1. The highest weight $\chi$ of the representation induced from
$\p_{\alpha}$ can be uniquely determined with an accuracy to within
$\chi (1)$) from $\chi (h)$.

2. For $\alpha^2=\lambda (\lambda +2)$ the corresponding subalgebra
$\p_{\alpha}\subset {\frak s}$ is not maximal: it is contained
in the subalgebra
$\p'_{\alpha}=\n_+\oplus \h\oplus \bigl\{ f^k(h-\alpha )p(h)$,
$p(h)\in  \C \ [h],\ k\geq 1\bigr\}$.
\endremark

The remaining part of Sec.~3 is devoted to the proof of the theorem.
\subsubhead
3.2
\endsubsubhead
Theorem~1 is a simple corollary of the following weaker result.

\proclaim{Theorem~2}
For a general $\alpha$, $\chi (h)$ the representatiuon of the Lie algebra
${\frak s}$ induced from $\p_{\alpha}$ with the highest weight
$\chi \in \h^*$ is irreducible.
\endproclaim

We saw in Sec.~2 that the irreducibility for the general
$\alpha$, $\chi (h)$ is equivalent to the global identity which can be
violated only if conditions~(5), (i)--(ii) are satisfied by a not
maximal of points (see Subsec.~2.4), i.e., if $T_{1,t}\bigl(
{\displaystyle\frac{\alpha}{t}}+2\bigr) ={\displaystyle\frac{1}{2}}\bigl(
-{\displaystyle\frac{\alpha}{t}}-{\displaystyle\frac{\alpha^2}{2t^2}}+{
\displaystyle\frac{\lambda (\lambda +2)}{2t^2}}\bigr)$ does not contain
terms with ${\displaystyle\frac{1}{t^2}}$, i.e.,
$\alpha^2=\lambda (\lambda +2)$,
or for $\chi (h)=0$. It is clear that the nonsatisfaction of this identity
for  certain $\alpha$, $\chi (h)$ means that the determinant of Shapovalov's
form of the corresponding representation vanishes. Therefore Theorem~1
follows from Theorem~2.

Below we shall prove Theorem~2.
\subsubhead
3.3
\endsubsubhead
Since the determinant of Shapovalov's form is an analytic function of
$\lambda \in \C$, Theorem~2 is a corollary of the following theorem
on the representations of the Lie algebra ${\frak s}_0$.

\proclaim{Theorem~3} The representations of the Lie algebra ${\frak s}_0$,
induced from $\p_{\alpha}$, with the highest weight $\chi \in \h^*$
are irreducible for  $\alpha$, $\chi (h)\ne 0$.
\endproclaim

The proof of Theorem~3 is given in Subsecs.~3.4--3.5.
\subsubhead
3.4
\endsubsubhead
The functions on the cone $\bigl\{ 2ef
+{\displaystyle\frac{h^2}{2}}=0\bigr\}$, which have the point $0\in \C^{\
3}$ as a zero of order $\geq k$, form an ideal
$I_k\subset {\frak s}_0$. The Lie algebra $I_2$ is nilpotent, and
$I_3=\{ I_2,I_2\} ,\ I_4=\{ I_2,I_3\}$ etc.

For nilpotent finite-dimensional Lie algebra Kirillov's theorem~[5]
states that the representations induced from the maximal parabolic
subalgebras are irreducible. To be more precise, let $\G$ be such a
Lie algebra, $\chi \in \G^*$. Let ${\frak b}$ be a maximal
subalgebra in $\G$ on which the bilinear form $\langle x,y\rangle =\chi
(\{ x,y\} ),\ x,y\in \G$ is zero. Then $\chi$ defines the one-dimensional
${\frak b}$-module
$\C_{\ \chi}$ and the induced representation $U(\G
)\bigotimes\limits_{U({\frak b} )}\C_{\ \chi}$ is irreducible.
We shall prove this theorem in Appendix~A.

Suppose now that $\chi \in \h^*$ is the hifhest weight of the
representation of the algebra ${\frak s}_0$ induced from $\p_{\alpha}$.
We continue it to the functional on ${\frak s}_0$ by setting
$\chi \vert_{\n_+}=\chi \vert_{\n_-}=0$.

\proclaim{Lemma~2} Let $\alpha$, $\chi (h)\ne 0$. Then the intersection
$I_2\cap \p_{\alpha}$ is the maximal subalgebra in $I_2$ on which the form
$\langle x,y\rangle =\chi \vert_{I_2}(\{ x,y\} );  x,y\in
I_2$ is zero.
\endproclaim

\demo{Proof}
It is obvious that
$I_2\cap \p_{\alpha}=(I_2\cap \n_+)\oplus (I_2\cap \h)\oplus
\{ f^k(h-\alpha )p(h)$, $k\geq 2$, $p(h)\in \C \ [h]\}
\oplus \{ fh(h-\alpha )p(h),\ p(h)\in \C \ [h]\}$. Let
${\widetilde \p}_{\alpha}$ be the maximal subalgebra in $I_2$ which
contains $\p_{\alpha}\cap I_2$ and on which the form $\langle x,y\rangle$
is zero.
We must prove that ${\widetilde \p}_{\alpha}=\p_{\alpha}$.

(i) Let $f^k(h-\alpha )^lp(h)\in {\widetilde \p}_{\alpha}$, where $l<k$ and
$p(h)$ is relatively prime to $(h-\alpha )$. Commuting this element with
$h^2,h^3,\ldots \in I_2\cap \p_{\alpha}$, we see that
$f^k(h-\alpha )^lp(h)p_1(h)\in {\widetilde \p}_{\alpha}$ for any
$p_1(h)\in \C \ [h]$.
In particular, $f^k(h-\alpha )^{k-l}\in {\widetilde
\p}_{\alpha}$.

(ii) If $fh\in {\widetilde \p}_{\alpha}$, then $\chi (\{ fh,e\} ) =0$
$\chi (\{ f(h-\alpha ),e\} )$ is also zero.
Hence contrary to the hypothesis of the lemma, either $\alpha$ or
$\chi (h)$ is equal to zero.

(iii) Assume that $f^k(h-\alpha )^{k-l}\in {\widetilde \p}_{\alpha}$ (see (i)),
$\{ f^k(h-\alpha )^{k-1},eh\} =
-kf^{k-1}h^2(h-\alpha)^{k-1}+(k-1)f^k(h-
\alpha )^{k-2}\cdot h\cdot e+2kf^ke(h-\alpha )^{k-1}$
lies in ${\widetilde \p}_{\alpha}$ since $eh\in \p_{\alpha}$.
The first and third terms lie in $\p_{\alpha}$ and, hence,
$f^{k-1}h^3(h-\alpha)^{k-2}$ lies in
${\widetilde \p}_{\alpha}$. According to (i), it
follows that $f^{k-1}(h-\alpha )^{k-2}\in {\widetilde \p}_{\alpha}$.
Thus, lowering the degree step by step,
we arrive at a contradiction with~(ii).
\enddemo
\subsubhead
{3.5. Proof of Theorem~2:}
\endsubsubhead

\proclaim{Lemma~3}
Let $\alpha$, $\chi (h)\ne 0$. Then the representation
$U(I_2)\bigotimes\limits_{U(I_2\cap
\p_{\alpha})}\C_{\ \chi \vert_{I_2}}$ of the Lie algebra
$I_2$ is irreducible.
\endproclaim

\demo{Proof}
In order to use Kirillov's theorem (see Subsec.~3.4),
it suffices to reduce the consideration to a finite-dimensional case.
The quotient algebra $I_2/I_k$
is nilpotent and finite-dimensional and the corresponding induced
representation of this Lie algebra is irreducible according to Kirillov's
theorem. However, for $k\gg 0$ the factorization does not affect the
levels which  are considerably smaller than $k$. In particular, there are
no $I_2$-singular vectors on them.
\enddemo
\remark{Remark}
Let ${\frak b}_k$ be the correspondinhg maximal subalgebra in
$I_2/I_k$. It is not true that
${\frak b}_k=I_2\cap \p_{\alpha}\big/ I_k\cap
\p_{\alpha}$; in particular,
$I_{k-1}\big/ I_k\subset {\frak b}_k$. However, it is true for the
intersections with $\n_-^{(l)}$, where $l\ll k$.
\endremark

The space of the representation of the Lie algebra
${\frak s}_0$ induced from $\p_{\alpha}$ can be written as
$$
V_1 = \{ \text{monomials of}\ f^i h^j,\ j<i\},
$$
or as
$$
V_2 = \{ \text{monomials of}\ f^ih^j,\ i\geq 2,\ j<i,\quad \text{and}\ fh\} .
$$
There exists an isomorphism $V_1{\widetilde \to}V_2$ which commutes
with the action of ${\frak s}_0$. However, the restriction of
$V_2$ to $I_2$ coincides with the representation of ${\frak s}_0$, induced from
$I_2\cap \p_{\alpha}$, which is irreducible according to Lemma~3 if
$\alpha$, $\chi (h)\ne 0$.
\subhead
Appendix~A
\endsubhead
The following result due to A.~A.~Kirillov~[5]:
\proclaim
{Theorem}
Let $g$ be a fihite dimensional nilpotent Lie algebra $\chi\in g^*$, and let
$\frak b$ be a maximal subalgebra in $g$ such that the restriction of the
bilinear form $(x\mid y)=\chi\bigl([x,y]\bigr)$ on $g$ to $\frak b$ is zero.
Then the inducrd representation $M_{\chi}=U(g)\otimes_{U(\frak b)}
\C_{\chi}$ is irreducible.
\endproclaim

To make the reasoning of Chapter~3 self-contained, we shall prove it
in this Appendix.

We begint with with obvious lemmas.

\proclaim
{Lemma~1}
Denote ${\frak b}_0=\frak b$  and ${\frak b}=\text{norm}{\frak b}_{i-1}$.
Then there exists an integer $n$ such that ${\frak b}_n=g$ and inclusions
$\frak b=\frak b\subset{\frak b}_1\subset\ldots\subset{\frak b}_n=g$
are strict.
\endproclaim
\proclaim{Lemma~2}
If $x\in{\frak b}_i$, $y={\frak b}_j$ $(i,\,j\ge 1)$ then
$[x,y]\in{\frak b}_{i+j-1}$.
\endproclaim

Using filtration on $g$ from Lemma~1, we can define a filtration
$$
\C=U_0\subset U_1\subset U_2\subset\ldots
$$
on the universal enveloping algebra $U(g)$. It follows easily from Lemma~2
that the ajoint algebra corresponding to this filtration is commutative.
Hence, any $b\in\frak b$ defines the mapping
$$
\operatorname{ad}(b):U_k/U_{k-1}\to U_{k-1}/U_{k-2}.
$$

\proclaim{Lemma~3}
For any $u\in U_k/U_{k-1}$, $k\ge 2$, there exists $b\in \frak b$ such that
$\operatorname{ad}(b)(u)\ne 0$ in $U_{k-1}/U_{k-2}$.
\endproclaim

\demo{Proof of the Theorem:}

Let  $v$ be the highest vector in $M_{\chi}$. Supposr that for some
$x\in{\frak b}_1$ $xv$ lies in a proper submodule of $M_{\chi}$. Then there
exists $b\in\frak b$ for which $\chi\bigl([x,b]\bigr)\ne 0$. Indeed, if
$\chi\bigl([x,b]\bigr)=0$ for all $b\in\frak b$ we can join $x$ to $\frak b$
and, hence, $\frak b$ is not maximal.

Then $bv=\alpha v$ for some $\alpha\in \C$, and $[b,x]v=bxv-\alpha x v$ lies in
a proper submodule and, hence, this submodule consides with $M_{\chi}$,
contrary to the fact that it is proper.

Now, assume that there exists a proper submodule $M$ in $M_{\chi}$, and let
$u\cdot v$ be some element of $M$, where $u\in U_k\setminus U_{k-1}$.
According to Lemma~3, there exists $b\in\frak b$ such that
$\operatorname{ad(b)(u)\ne 0}$ as an element of $U_{k-1}/U_{k-2}$, and
therefore $\operatorname{ad}(b)(u)v$ is a nonzero element of $M$ with
lower filtration. Iterating this process, we shall arrive to the case
in the beginning of the proof.
\enddemo
\subhead
Appendix~B
\endsubhead
In fact, using methods of Ch.~1 one can easily
prove the following generalization of identity~(16) of Ch.~1
and of the local identity (see 0.11)
$$
\prod_{i\ge 1}\frac1{(1-a\cdot q^i)^i}
=\sum\Sb{\text{over all}}\\{\text{Young diagrams}}\\
{\Cal D}_{l_1,\ldots,l_k}\endSb a^{\#{\Cal D}_{l_1,\ldots,l_k}}\cdot
q^{\sum l_i\cdot i^2}\cdot\bigl(\chi({\Cal
D}_{l_1,\ldots,l_k})\bigr)^2
\tag 1
$$
(the notation is as in 0.11). Therefore, in the reasoning of
Sec.~4 of Chapter~2 imbedding into
$gl_\infty\bigl(\C[t]/t^2\bigr)$
can be replaced by more elementary considerations, using
additional grading.
This proves both the local identity and its  generalization~(1)
for higher derivators. For the time  being, the author does not
know any generalizations of the global identity (see 0.11)
for higher derivators.

\enddocument